\documentclass[12pt]{article}
\usepackage[]{amsmath,amssymb}
\usepackage{amscd}
\usepackage{latexsym}
\usepackage{cite}

\newtheorem{definition}{Definition}[section]
\newtheorem{theorem}[definition]{Theorem}
\newtheorem{lemma}[definition]{Lemma}
\newtheorem{corollary}[definition]{Corollary}

\newtheorem{example}[definition]{Example}
\newtheorem{conjecture}[definition]{Conjecture}
\newtheorem{problem}[definition]{Problem}
\newtheorem{note}[definition]{Note}

\newtheorem{proposition}[definition]{Proposition}

\def\F{\mathbb F}
\def\K{\mathbb F}
\def\Knew{\mathbb K}
\def\vare{\varepsilon}

\begin{document}
\title{\bf A classification of sharp \\tridiagonal pairs}
\author{
Tatsuro Ito\footnote{Supported in part by JSPS grant
18340022.}, 
Kazumasa Nomura,
and
Paul Terwilliger}
\date{}

\maketitle
\begin{abstract}
\noindent
Let $\F$ denote a 
 field and let $V$ denote a vector space over $\F$ with
 finite positive dimension.
 We consider a pair of linear transformations $A:V \to V$
 and $A^*:V \to V$ that satisfy the following conditions:
 (i)
 each of $A,A^*$ is diagonalizable;
 (ii)
 there exists an ordering $\lbrace V_i\rbrace_{i=0}^d$ of the eigenspaces of
 $A$ such that
 $A^* V_i \subseteq V_{i-1} + V_{i} + V_{i+1}$ for $0 \leq i \leq d$,
 where $V_{-1}=0$ and $V_{d+1}=0$;
 (iii)
 there exists an ordering $\lbrace V^*_i\rbrace_{i=0}^\delta$
 of the eigenspaces of $A^*$ such that
 $A V^*_i \subseteq V^*_{i-1} + V^*_{i} + V^*_{i+1}$ for
  $0 \leq i \leq \delta$,
  where $V^*_{-1}=0$ and $V^*_{\delta+1}=0$;
  (iv)
  there is no subspace $W$ of $V$ such that
  $AW \subseteq W$, $A^* W \subseteq W$, $W \neq 0$, $W \neq V$.
  We call such a pair a {\it tridiagonal pair} on $V$.
  It is known that $d=\delta$ and for
  $ 0 \leq i \leq d$ the dimensions of
  $V_i,V_{d-i},V^*_i, V^*_{d-i}$ coincide. The pair
  $A,A^*$ is called {\it sharp} whenever ${\rm dim} \,V_0=1$.
It is known that if $\F$ is algebraically closed then
$A,A^*$ is sharp.
In this paper we classify up to isomorphism the
sharp tridiagonal pairs. As a corollary,
we classify  up to isomorphism the tridiagonal pairs
over an algebraically closed field.
We obtain these classifications by proving
 the $\mu$-conjecture.

 \end{abstract}

\section{Tridiagonal pairs}

\noindent 
Throughout this paper $\K$ denotes a field
and $\overline \K$ denotes the algebraic closure of $\K$.
An algebra is meant to be associative and have a $1$.

\medskip
\noindent 
We begin by recalling the notion of a tridiagonal pair. 
We will use the following terms.
Let $V$ denote a vector space over $\K$ with finite
positive dimension.
For a 
 linear transformation $A:V\to V$
and a
subspace $W \subseteq V$,
we call $W$ an
 {\it eigenspace} of $A$ whenever 
 $W\not=0$ and there exists $\theta \in \K$ such that 
$W=\lbrace v \in V \;\vert \;Av = \theta v\rbrace$;
in this case $\theta$ is the {\it eigenvalue} of
$A$ associated with $W$.
We say that $A$ is {\it diagonalizable} whenever
$V$ is spanned by the eigenspaces of $A$.

\begin{definition}  
{\rm \cite[Definition~1.1]{TD00}}
\label{def:tdp}
\rm
Let $V$ denote a vector space over $\K$ with finite
positive dimension. 
By a {\it tridiagonal pair} (or {\it $TD$ pair})
on $V$
we mean an ordered pair of linear transformations
$A:V \to V$ and 
$A^*:V \to V$ 
that satisfy the following four conditions.
\begin{enumerate}
\item Each of $A,A^*$ is diagonalizable.
\item There exists an ordering $\lbrace V_i\rbrace_{i=0}^d$ of the  
eigenspaces of $A$ such that 
\begin{equation}
A^* V_i \subseteq V_{i-1} + V_i+ V_{i+1} \qquad \qquad 0 \leq i \leq d,
\label{eq:t1}
\end{equation}
where $V_{-1} = 0$ and $V_{d+1}= 0$.
\item There exists an ordering $\lbrace V^*_i\rbrace_{i=0}^{\delta}$ of
the  
eigenspaces of $A^*$ such that 
\begin{equation}
A V^*_i \subseteq V^*_{i-1} + V^*_i+ V^*_{i+1} 
\qquad \qquad 0 \leq i \leq \delta,
\label{eq:t2}
\end{equation}
where $V^*_{-1} = 0$ and $V^*_{\delta+1}= 0$.
\item There does not exist a subspace $W$ of $V$ such  that $AW\subseteq W$,
$A^*W\subseteq W$, $W\not=0$, $W\not=V$.
\end{enumerate}
We say the pair $A,A^*$ is {\it over $\K$}.
We call $V$ the 
{\it underlying
 vector space}.
\end{definition}

\begin{note}
\label{lem:convention}
\rm
According to a common notational convention $A^*$ denotes 
the conjugate-transpose of $A$. We are not using this convention.
In a TD pair $A,A^*$ the linear transformations $A$ and $A^*$
are arbitrary subject to (i)--(iv) above.
\end{note}

\medskip
\noindent 
We now give some background on TD pairs;
for more information we refer the reader to the
survey
\cite{madrid}.
The concept of a TD pair
originated in algebraic graph theory,
or more precisely, the theory of
$Q$-polynomial distance-regular graphs.
The concept is implicit in
\cite[p.~263]{BanIto},
\cite{Leon}
and more explicit in 
\cite[Theorem~2.1]{TersubI}.
A systematic study began in
\cite{TD00}.
Some notable papers on the topic are
\cite{
atakdual,
bas6,
shape,
tdanduq,
NN,
qtet,
IT:aug,
drin,
IT:qRacah,
LS99
}.
There are connections to 
representation theory
\cite{
bas1,hasan2,
neubauer,
Ha,
tdanduq,
qtet,
Koor,
Koor2,
Rosengren,
Rosengren2,
IT:Mock,
aw},
partially ordered sets
    \cite{lsint},
the bispectral problem 
\cite{atak, 
atakdual,
GYLZmut, GH7,GH1,Zhidd},
statistical mechanical models
 \cite{bas1,bas2,bas3,bas4,bas5,bas6,bas7,DateRoan2, Davfirst, 
Dolgra,
 Da,Onsager},
and other areas of physics
\cite{odake, VZ, LPcm}.

\medskip
\noindent 
Let $A,A^*$ denote a TD pair
on $V$, as in Definition 
\ref{def:tdp}. By
\cite[Lemma 4.5]{TD00}
the integers $d$ and $\delta$ from
(ii), (iii) are equal; we call this
common value the {\it diameter} of the
pair.
By \cite[Theorem 10.1]{TD00} the pair $A, A^*$ satisfy two polynomial
equations called the tridiagonal relations;
these generalize the $q$-Serre relations
  \cite[Example~3.6]{qSerre}
and the Dolan-Grady relations
  \cite[Example~3.2]{qSerre}.
See
\cite{
bas6,
GH1,
IT:aug,
IT:qRacah,
N:aw,
tersub3,
qSerre,
aw,
vidunas,
V2}
 for results on the tridiagonal relations.
An ordering of the eigenspaces of $A$ (resp. $A^*$)
is said to be {\em standard} whenever it satisfies 
(\ref{eq:t1})
 (resp. (\ref{eq:t2})). 
We comment on the uniqueness of the standard ordering.
Let $\{V_i\}_{i=0}^d$ denote a standard ordering of the eigenspaces of $A$.
By \cite[Lemma~2.4]{TD00}, 
 the ordering $\{V_{d-i}\}_{i=0}^d$ is also standard and no further
 ordering
is standard.
A similar result holds for the eigenspaces of $A^*$.
Let $\{V_i\}_{i=0}^d$ (resp.
$\{V^*_i\}_{i=0}^d$)
denote a standard ordering of the eigenspaces
 of $A$ (resp. $A^*$).
For $0 \leq i \leq d$ let 
$\theta_i$ 
(resp. $\theta^*_i$)
denote the eigenvalue of
$A$ 
(resp.  $A^*$) associated with
$V_i$ 
(resp. $V^*_i$).
By
 \cite[Theorem~11.1]{TD00}
 the expressions
\begin{eqnarray*}
\frac{\theta_{i-2}-\theta_{i+1}}{\theta_{i-1}-\theta_i},  \qquad\qquad
  \frac{\theta^*_{i-2}-\theta^*_{i+1}}{\theta^*_{i-1}-\theta^*_i}
\end{eqnarray*}
are equal and independent of $i$ for $2 \leq i \leq d-1$.
We call the
sequence
$\lbrace \theta_i\rbrace_{i=0}^d$
(resp. $\lbrace \theta^*_i\rbrace_{i=0}^d$)
the
{\it eigenvalue sequence}
(resp. {\it  dual eigenvalue sequence})
for  the given  standard orderings.
See 
\cite{TD00,drin,Pas, 
LS99,
  qSerre} for results on the eigenvalues and dual eigenvalues.
By \cite[Corollary~5.7]{TD00}, 
for $0 \leq i \leq d$ the spaces $V_i$, $V^*_i$
have the same dimension; we denote
this common dimension by $\rho_i$. 
By \cite[Corollaries 5.7, 6.6]{TD00}
the sequence $\{\rho_i\}_{i=0}^d$ is symmetric    and unimodal;
that is $\rho_i=\rho_{d-i}$ for $0 \leq i \leq d$ and
$\rho_{i-1} \leq \rho_i$ for $1 \leq i \leq d/2$.
By \cite[Theorem~1.3]{NT:onshape} 
we have
$\rho_{i} \leq 
\rho_{0}\binom{d}{i}
$ for $0 \leq i \leq d$.
We call the sequence $\{\rho_i\}_{i=0}^d$ the {\em shape}
of $A,A^*$.
See 
\cite{
shape,
sharpen,
nom3,
nom4,
nomsharp,
NT:onshape}
for results on the shape.
The TD pair $A,A^*$ is called {\it sharp} whenever
$\rho_0=1$.
By
\cite[Theorem~1.3]{nomstructure},
if $\F$ is algebraically closed then
$A,A^*$ is sharp.
In any case $A,A^*$ can be
``sharpened'' by replacing $\F$ 
with a certain field extension $\Knew$
of $\F$ 
that has index
$\lbrack \Knew:\F \rbrack=\rho_0$
\cite[Theorem~4.12]{sharpen}.
Suppose that $A,A^*$ is sharp.
Then by \cite[Theorem~1.4]{nomstructure}, there exists a
nondegenerate symmetric bilinear form $\langle \,,\,\rangle $
on $V$ such that
$\langle Au,v \rangle =\langle u,Av \rangle$
and $\langle A^*u,v \rangle =\langle u,A^*v \rangle $
for all $u,v \in V$.
See \cite{CurtH,nomsharp,Tanaka2,Tang} for results on the bilinear form.

\medskip
\noindent 
The following special cases of TD pairs have been studied extensively.
In \cite{Vidar} the TD pairs of shape $(1,2,1)$ are classified and
described in detail.
A TD pair of shape
$(1,1,\ldots, 1)$
is called a {\it Leonard pair}
 \cite[Definition 1.1]{LS99}, and these 
are classified 
in  \cite[Theorem 1.9]{LS99}.
This classification yields a  
 correspondence between the Leonard pairs and a
family of orthogonal polynomials consisting of the $q$-Racah polynomials
and their relatives 
\cite{AWil,
TLT:array,
qrac}.
This family coincides with the terminating branch of the Askey scheme
\cite{KoeSwa}.
See 
\cite{Cur1, Cur2, Mik2, NT:balanced,NT:formula,NT:det,NT:mu,
NT:span,NT:switch,
nomsplit,
madrid}
and the references therein for 
results on Leonard pairs.
Our TD  pair $A,A^*$ is said to have {\em Krawtchouk type}
(resp. {\em $q$-geometric type})
whenever $\lbrace d-2i \rbrace_{i=0}^d $
(resp. $\lbrace q^{d-2i}\rbrace_{i=0}^d$)
is both an eigenvalue sequence and dual eigenvalue sequence for
the pair.
In \cite[Theorems 1.7, 1.8]{Ha} Hartwig classified the
TD pairs over $\F$ that have Krawtchouk type, provided that
$\F$ is algebraically closed with characteristic zero. 
By \cite[Remark~1.9]{Ha} these TD pairs  are in bijection with the
finite-dimensional irreducible modules for the  three-point loop algebra
$\mathfrak{sl}_2 \otimes \K[t,t^{-1}, (t-1)^{-1}]$.
See
\cite{
Ha,
HT,
Ev,
IT:Krawt,
drin
}
for results on  TD pairs of Krawtchouk type.
In \cite[Theorems~1.6,~1.7]{NN} we classified 
the TD pairs over $\F$ that have $q$-geometric type,
provided that
$\F$ is algebraically closed and $q$ is not a root of unity.
By \cite[Theorems~10.3,~10.4]{qtet} these TD pairs are in bijection
with the type $1$, finite-dimensional, irreducible modules for the 
$\F$-algebra
$\boxtimes_q$; this is a $q$-deformation of
$\mathfrak{sl}_2 \otimes \K[t,t^{-1}, (t-1)^{-1}]$ as explained in \cite{qtet}.
See
\cite{hasan,
hasan2,
shape,
tdanduq,
NN,
qtet,
ITdrg,
drin
} 
for results  on  $q$-geometric
TD pairs. 
There is a general family of TD pairs said to have  $q$-Racah type;
these have an eigenvalue sequence 
and dual eigenvalue sequence of the form
(\ref{eq:const1})--(\ref{eq:const5}) below.
The Leonard pairs of $q$-Racah type correspond to the
$q$-Racah polynomials 
\cite[Example~5.3]{TLT:array}.
In \cite[Theorem~3.3]{IT:qRacah} we classified 
the TD pairs over $\F$ that have $q$-Racah type,
provided that $\F$ is algebraically closed.
See 
\cite{IT:aug, drin, IT:qRacah} for results 
on 
TD pairs of $q$-Racah type.

\medskip
\noindent Turning to the
present paper, in our main result we classify up to isomorphism
the sharp TD pairs.
Here is a summary of the argument.
In 
\cite[Conjecture~14.6]{IT:Krawt}
we conjectured how a
classification of all the sharp TD pairs
would look;
this is the
{\it classification conjecture}. 
Shortly afterwards we introduced a conjecture,
called the {\em $\mu$-conjecture},
which implies the classification conjecture. The $\mu$-conjecture is roughly
described as follows. Start with a sequence
$p=(\{\theta_i\}_{i=0}^d;\{\theta^*_i\}_{i=0}^d)$
of scalars taken from $\K$ that
satisfy the known constraints on the eigenvalues of a TD pair over
$\K$ of diameter $d$; these are
conditions (i), (ii) in Theorem
\ref{conj:main} below.
Following \cite[Definition 2.4]{nomstructure} we associate with $p$
an  $\K$-algebra $T$ defined by generators
and relations; see Definition
\ref{def:twitha} for the precise definition.
We are interested in the $\K$-algebra $e^*_0Te^*_0$ where $e^*_0$ is a certain
idempotent element of $T$. Let $\{x_i\}_{i=1}^d$ denote mutually commuting
indeterminates. Let $\K[x_1,\ldots,x_d]$ denote the $\K$-algebra consisting of
the polynomials in $\{x_i\}_{i=1}^d$ that have all coefficients in $\K$.
In \cite[Corollary 6.3]{nom:mu} we displayed a surjective
$\K$-algebra homomorphism $\mu: \K[x_1,\ldots,x_d] \to e^*_0Te^*_0$.
The {\em $\mu$-conjecture} \cite[Conjecture 6.4]{nom:mu} asserts that
$\mu$ is an isomorphism.
By \cite[Theorem 10.1]{nom:mu}
the $\mu$-conjecture implies the classification
conjecture.
In 
\cite[Theorem 12.1]{nom:mu} 
we showed that the $\mu$-conjecture holds for $d \leq 5$.
In 
\cite[Theorem~5.3]
{NT:muqrac}
we showed that the
$\mu$-conjecture holds for the case
in which $p$ has $q$-Racah type.
In the present paper we combine this fact with some algebraic geometry
to prove the $\mu$-conjecture in general. The
 $\mu$-conjecture (now a theorem) is given in
Theorem
\ref{conj:mainp}.
Theorem 
\ref{conj:mainp} 
implies the
classification conjecture, and this yields our
classification of the sharp TD pairs.
The classification is given in Theorem 
 \ref{conj:main}.
As a corollary, we classify up to isomorphism the
TD pairs over an algebraically closed field.
This result can be found in
Corollary
 \ref{conj:mainac}.


\medskip
\noindent 
Section 3 contains
the precise statements
of our main results.
 In Section 2 we review the concepts needed to
make these statements.

\section{Tridiagonal systems}

\indent
When working with a TD pair, it is often convenient to consider
a closely related object called a TD system.
To define a TD system, we recall a few concepts from linear
algebra.
Let $V$ denote a vector space over $\K$ with finite
positive dimension.
Let ${\rm End}(V)$ denote the $\K$-algebra of all linear
transformations from $V$ to $V$.
Let $A$ denote a diagonalizable element of $\mbox{\rm End}(V)$.
Let $\{V_i\}_{i=0}^d$ denote an ordering of the eigenspaces of $A$
and let $\{\theta_i\}_{i=0}^d$ denote the corresponding ordering of
the eigenvalues of $A$.
For $0 \leq i \leq d$ define $E_i \in 
\mbox{\rm End}(V)$ 
such that $(E_i-I)V_i=0$ and $E_iV_j=0$ for $j \neq i$ $(0 \leq j \leq d)$.
Here $I$ denotes the identity of $\mbox{\rm End}(V)$.
We call $E_i$ the {\em primitive idempotent} of $A$ corresponding to $V_i$
(or $\theta_i$).
Observe that
(i) $I=\sum_{i=0}^d E_i$;
(ii) $E_iE_j=\delta_{i,j}E_i$ $(0 \leq i,j \leq d)$;
(iii) $V_i=E_iV$ $(0 \leq i \leq d)$;
(iv) $A=\sum_{i=0}^d \theta_i E_i$.
Here $\delta_{i,j}$ denotes the Kronecker delta.
Note that
\begin{equation}         \label{eq:defEi}
  E_i=\prod_{\stackrel{0 \leq j \leq d}{j \neq i}}
          \frac{A-\theta_jI}{\theta_i-\theta_j}
\qquad \qquad 0 \leq i \leq d.
\end{equation}
Observe that each of 
$\{A^i\}_{i=0}^d$,
$\{E_i\}_{i=0}^d$ is a basis for the $\K$-subalgebra
of $\mbox{\rm End}(V)$ generated by $A$.
Moreover $\prod_{i=0}^d(A-\theta_iI)=0$.
Now let $A,A^*$ denote a TD pair on $V$.
An ordering of the primitive idempotents 
 of $A$ (resp. $A^*$)
is said to be {\em standard} whenever
the corresponding ordering of the eigenspaces of $A$ (resp. $A^*$)
is standard.

\begin{definition}
{\rm \cite[Definition~2.1]{TD00}}
 \label{def:TDsystem} 
\rm
Let $V$ denote a vector space over $\K$ with finite
positive dimension.
By a {\it tridiagonal system} (or {\it  $TD$ system}) on $V$ we mean a sequence
\[
 \Phi=(A;\{E_i\}_{i=0}^d;A^*;\{E^*_i\}_{i=0}^d)
\]
that satisfies (i)--(iii) below.
\begin{itemize}
\item[(i)]
$A,A^*$ is a TD pair on $V$.
\item[(ii)]
$\{E_i\}_{i=0}^d$ is a standard ordering
of the primitive idempotents of $A$.
\item[(iii)]
$\{E^*_i\}_{i=0}^d$ is a standard ordering
of the primitive idempotents of $A^*$.
\end{itemize}
We say that $\Phi$ is {\em over} $\K$.
We call $V$ the {\it underlying vector space}.
\end{definition}

\noindent The following result is immediate 
from
lines
(\ref{eq:t1}), 
(\ref{eq:t2}) 
and Definition
 \ref{def:TDsystem}. 

\begin{lemma}
\label{lem:triplep}
Let $(A; \lbrace E_i\rbrace_{i=0}^d; A^*; \lbrace E^*_i\rbrace_{i=0}^d)$
denote a TD system. Then the following hold for $0 \leq i,j,k\leq d$.
\begin{enumerate}
\item[\rm (i)] $E^*_iA^kE^*_j=0$ if $k < |i-j|$;
\item[\rm (ii)] $E_iA^{*k}E_j=0$ if $k < |i-j|$.
\end{enumerate}
\end{lemma}

\medskip
\noindent The notion of isomorphism for TD systems
is defined in
\cite[Section~3]{nomsharp}.

\begin{definition}        \label{def}
\rm
Let $\Phi=(A;\{E_i\}_{i=0}^d;A^*$; $\{E^*_i\}_{i=0}^d)$ 
denote a TD system on $V$.
For $0 \leq i \leq d$ let $\theta_i$ (resp. $\theta^*_i$)
denote the eigenvalue of $A$ (resp. $A^*$)
associated with the eigenspace $E_iV$ (resp. $E^*_iV$).
We call $\{\theta_i\}_{i=0}^d$ (resp. $\{\theta^*_i\}_{i=0}^d$)
the {\em eigenvalue sequence}
(resp. {\em dual eigenvalue sequence}) of $\Phi$.
Observe that $\{\theta_i\}_{i=0}^d$ (resp. $\{\theta^*_i\}_{i=0}^d$) are
mutually distinct and contained in $\K$.
We call $\Phi$ 
{\it sharp} whenever the TD pair
$A,A^*$ is sharp.
\end{definition}
\medskip


\noindent 
The following notation will be useful.

\begin{definition}
\rm
\label{def:taueta}
Let $x$ denote an indeterminate and let $\K[x]$
denote the $\K$-algebra consisting of the polynomials
in $x$ that have all coefficients in $\K$.
Let $\lbrace \theta_i\rbrace_{i=0}^d$  and
 $\lbrace \theta^*_i\rbrace_{i=0}^d$
denote scalars in $\K$.
For $0 \leq i \leq d$  define the following polynomials in
$\K[x]$:
\begin{eqnarray*}
 \tau_i &=& 
  (x-\theta_0)(x-\theta_1)\cdots(x -\theta_{i-1}), \\
 \eta_i &=&
  (x-\theta_d)(x-\theta_{d-1})\cdots(x-\theta_{d-i+1}),  \\
 \tau^*_i &=&
  (x-\theta^*_0)(x-\theta^*_1)\cdots(x-\theta^*_{i-1}), \\
 \eta^*_i &=&
  (x-\theta^*_d)(x-\theta^*_{d-1})\cdots(x-\theta^*_{d-i+1}).
\end{eqnarray*}
Note that each of $\tau_i$, $\eta_i$, $\tau^*_i$, $\eta^*_i$ is
monic with degree $i$.
\end{definition}

\noindent 
We now recall the split sequence of a sharp TD system.
This sequence was originally defined
in \cite[Section~5]{IT:Krawt} using the split
decomposition \cite[Section~4]{TD00},
 but
in \cite{nom:mu}
an alternate
definition was introduced that is more convenient
to our purpose.

\begin{definition}
\label{def:split}
\rm
{\rm \cite[Definition~2.5]{nom:mu}}
Let 
$(A; \lbrace E_i\rbrace_{i=0}^d; A^*; \lbrace E^*_i\rbrace_{i=0}^d)$
denote a sharp TD system over $\K$, with eigenvalue
sequence $\lbrace \theta_i \rbrace_{i=0}^d$
and dual eigenvalue sequence
 $\lbrace \theta^*_i \rbrace_{i=0}^d$.
By 
\cite[Lemma 5.4]{nomstructure},
 for $0 \leq i \leq d$
there exists a unique $\zeta_i \in \K$ such that 
\begin{eqnarray*}
E^*_0 \tau_i(A) E^*_0 = 
\frac{\zeta_i E^*_0}{
(\theta^*_0-\theta^*_1) 
(\theta^*_0-\theta^*_2) 
\cdots
(\theta^*_0-\theta^*_i)}. 
\end{eqnarray*}
Note that $\zeta_0=1$.
We call $\lbrace \zeta_i \rbrace_{i=0}^d$
the {\it split sequence} of the TD system.
\end{definition}

\begin{definition}
\label{def:pa}
\rm
Let $\Phi$
denote a sharp TD system.
By the {\it parameter array} of $\Phi$ we mean the
sequence
 $(\{\theta_i\}_{i=0}^d; \{\theta^*_i\}_{i=0}^d; \{\zeta_i\}_{i=0}^d)$
where
 $\{\theta_i\}_{i=0}^d$
(resp. 
$\{\theta^*_i\}_{i=0}^d$
)
is the eigenvalue sequence 
(resp. dual eigenvalue sequence)
of $\Phi$ and
$\{\zeta_i\}_{i=0}^d$ is the split sequence of $\Phi$.
\end{definition}

\noindent 
The following result shows the significance
of the parameter array. 

\begin{proposition}
\label{thm:isopa}
{\rm 
\cite{IT:aug},
\cite[Theorem~1.6]{nomstructure}}
Two sharp TD systems over $\K$ are isomorphic
if and only if they have the same parameter array.
\end{proposition}

\section{Statement of results} 

\noindent In  this section we state our main results.
The first result below resolves 
\cite[Conjecture~14.6]{IT:Krawt}.

\begin{theorem}
 \label{conj:main}  
Let $d$ denote a nonnegative integer and  let
\begin{equation}         \label{eq:parrayc}
 (\{\theta_i\}_{i=0}^d; \{\theta^*_i\}_{i=0}^d; \{\zeta_i\}_{i=0}^d)
\end{equation}
denote a sequence of scalars taken from $\K$.
Then there exists a sharp TD system $\Phi$ over $\K$
with parameter array \eqref{eq:parrayc} if and only if
{\rm (i)--(iii)} hold below.
\begin{itemize}
\item[\rm (i)]
$\theta_i \neq \theta_j$, 
$\theta^*_i \neq \theta^*_j$ if $i \neq j$ $(0 \leq i,j \leq d)$.
\item[\rm (ii)]
The expressions
\begin{equation} 
\label{eq:betaplusone} 
\frac{\theta_{i-2}-\theta_{i+1}}{\theta_{i-1}-\theta_i},  \qquad\qquad
  \frac{\theta^*_{i-2}-\theta^*_{i+1}}{\theta^*_{i-1}-\theta^*_i}
\end{equation}
are equal and independent of $i$ for $2 \leq i \leq d-1$.
\item[\rm (iii)]
$\zeta_0=1$, $\zeta_d \neq 0$, and
\begin{equation*}   
        \label{eq:ineqc}
0 \neq \sum_{i=0}^d \eta_{d-i}(\theta_0)\eta^*_{d-i}(\theta^*_0) \zeta_i.
\end{equation*}
\end{itemize}
Suppose {\rm (i)--(iii)} hold. Then $\Phi$ is unique up to isomorphism of
TD systems.
\end{theorem}


\noindent In 
\cite[Conjecture~6.4]{nom:mu} we stated a conjecture
called the {\it $\mu$-conjecture}, and
we proved that the $\mu$-conjecture implies
Theorem
 \ref{conj:main}. To obtain 
 Theorem
 \ref{conj:main} we will prove the $\mu$-conjecture.
We now explain this conjecture.

\begin{definition}
\label{def:feas}
\rm
Let $d$ denote a nonnegative integer
and let 
 $(\lbrace \theta_i\rbrace_{i=0}^d;
\lbrace \theta^*_i\rbrace_{i=0}^d)$ 
denote a sequence of scalars taken from
$\F$.
This sequence is called {\it feasible} whenever
it satisfies
conditions
(i), (ii) of Theorem
 \ref{conj:main}. 
\end{definition}

\begin{definition}
\label{def:feaset}
\rm
For all integers $d\geq 0$ let 
${\rm Feas}(d,\F)$ 
denote the set of all
 feasible sequences
 $(\lbrace \theta_i\rbrace_{i=0}^d;
\lbrace \theta^*_i\rbrace_{i=0}^d)$ 
of scalars taken from $\F$.
\end{definition}

\begin{definition}
\label{def:twitha}
{\rm \cite[Definition~2.4]{nomstructure}}
\rm
Fix an integer $d\geq 0$ and
  a sequence $p=(\lbrace \theta_i\rbrace_{i=0}^d;
\lbrace \theta^*_i\rbrace_{i=0}^d)$  in 
${\rm Feas}(d,\F)$.
Let $T=T(p,\F)$ denote the $\F$-algebra defined by
generators
$a$, $\{e_i\}_{i=0}^d$, $a^*$, $\{e^*_i\}_{i=0}^d$ and relations
\begin{equation}                            \label{eq:eiejcop}
  e_ie_j=\delta_{i,j}e_i, \qquad
  e^*_ie^*_j=\delta_{i,j}e^*_i \qquad\qquad 
  0 \leq i,j \leq d,
\end{equation}
\begin{equation}               
\label{eq:sumeicop}
  1=\sum_{i=0}^d e_i,  \qquad\qquad
  1=\sum_{i=0}^d e^*_i, \quad
\end{equation}
\begin{equation}                                  \label{eq:sumthieicop}
   a = \sum_{i=0}^d \theta_ie_i, \qquad\qquad
   a^* = \sum_{i=0}^d \theta^*_i e^*_i,
\end{equation}
\begin{equation}                                     \label{eq:esiakesjcop}
 e^*_i a^k e^*_j = 0  \qquad \text{if $\;k<|i-j|$}
                  \qquad\qquad 0 \leq i,j,k \leq d,
\end{equation}
\begin{equation}                                      \label{eq:eiaskej}
 e_i {a^*}^k e_j = 0  \qquad \text{if $\;k<|i-j|$} 
                  \qquad\qquad 0 \leq i,j,k \leq d.
\end{equation}
\end{definition}

\begin{lemma}
{\rm \cite[Lemma~4.2]{nom:mu}}
In the algebra $T$ from Definition
\ref{def:twitha}, the elements $\lbrace e_i\rbrace_{i=0}^d$ are
linearly independent and the elements
 $\lbrace e^*_i\rbrace_{i=0}^d$ are
linearly independent.
\end{lemma}

\noindent The algebra $T$ is related to TD systems
as follows.
\begin{lemma}
{\rm \cite[Lemma~2.5]{nomstructure}}
Let $V$ denote a vector space over $\F$ with finite positive
dimension. Let
 $(A;\{E_i\}_{i=0}^d;A^*;\{E^*_i\}_{i=0}^d)$
 denote a TD system on $V$ with
 eigenvalue sequence 
 $\lbrace \theta_i\rbrace_{i=0}^d$ and
 dual eigenvalue sequence 
 $\lbrace \theta^*_i\rbrace_{i=0}^d$.
Let $T$ denote the $\K$-algebra from
Definition
\ref{def:twitha} corresponding to
 $(\lbrace \theta_i\rbrace_{i=0}^d;
 \lbrace \theta^*_i\rbrace_{i=0}^d)$.
Then there exists a unique $T$-module structure
on $V$ such that $a$, $e_i$, $a^*$, $e^*_i$ acts as
$A$, $E_i$, $A^*$, $E^*_i$ respectively. This
$T$-module is irreducible.
\end{lemma}

\noindent Fix an integer $d\geq 0$
 and a sequence 
$p \in {\rm Feas}(d,\F)$.
Let $T=T(p,\F)$ denote the corresponding algebra from
Definition
\ref{def:twitha}. Observe that
$e^*_0Te^*_0$ is an $\F$-algebra with multiplicative
identity $e^*_0$.

\begin{lemma} {\rm \cite[Theorem~2.6]{nomstructure}} 
\label{lem:17} 
With the above notation,
the algebra $e^*_0Te^*_0$ is commutative and generated by
\[
  e^*_0\tau_i(a)e^*_0 \qquad \qquad 1\leq i \leq d.
\]
\end{lemma}


\begin{corollary}  
{\rm \cite[Corollary~6.3]{nom:mu}}
\label{cor:18a}   
With the above notation,
there exists a surjective $\K$-algebra homomorphism 
$\mu: \K[x_1,\ldots,x_d]
\to e^*_0Te^*_0$ that sends $x_i \mapsto e^*_0\tau_i(a)e^*_0$
for $1 \leq i \leq d$.
\end{corollary}

\noindent In 
\cite[Conjecture~6.4]{nom:mu} we conjectured
that the map
$\mu$ from
Corollary
\ref{cor:18a}   
is an isomorphism. This is the
{\it $\mu$-conjecture}.
The following result resolves the $\mu$-conjecture.

\begin{theorem}    
\label{conj:mainp}  
Fix an integer $d\geq 0$
 and a  sequence
 $p \in
{\rm Feas}(d,\F)$. 
Let $T=T(p,\F)$ denote the corresponding algebra from
Definition
\ref{def:twitha}.
Then the map 
$\mu: \K[x_1,\ldots,x_d]
\to e^*_0Te^*_0$
 from Corollary  \ref{cor:18a} is an isomorphism.
\end{theorem}
\noindent We will prove 
Theorem 
 \ref{conj:main}  
and
Theorem \ref{conj:mainp}  
in Section 17.

\section{The $q$-Racah case}

\noindent Our proof of
Theorem \ref{conj:mainp}  
will use the fact that
the theorem is known to be true in a special case called
$q$-Racah
\cite[Theorem~5.3]{NT:muqrac}.
In this section we describe the
$q$-Racah case. We start with some comments
about the feasible sequences 
from Definition
\ref{def:feas}.


\begin{lemma}
\label{lem:pexist}
Assume $\F$ is infinite.
Then for all integers $d\geq 0$ the set 
${\rm Feas}(d,\F)$  is nonempty.
\end{lemma}
\noindent {\it Proof:}
Consider the
polynomial $\prod_{n=1}^d (x^n-1)$
in 
$\F\lbrack x \rbrack$.
Since
$\F$ is infinite there exists a nonzero $\vartheta \in \F$
that is not a root of this polynomial. Define
$\theta_i=\vartheta^i$ and 
$\theta^*_i=\vartheta^i$ for $0 \leq i \leq d$.
Then $(\lbrace \theta_i\rbrace_{i=0}^d;
\lbrace \theta^*_i\rbrace_{i=0}^d)$ 
satisfies the conditions (i), (ii) of
Theorem \ref{conj:main} and is therefore feasible.
The result follows.
\hfill $\Box$ \\

\noindent Fix an integer $d\geq 0$
and 
 a sequence 
$(\lbrace \theta_i\rbrace_{i=0}^d;
\lbrace \theta^*_i\rbrace_{i=0}^d)$ 
in 
${\rm Feas}(d,\F)$.
This sequence must satisfy 
condition 
(ii) in Theorem
\ref{conj:main}.
 For this constraint the ``most general'' solution is
\begin{eqnarray}
\label{eq:const1}
&&\theta_i = \alpha + b q^{2i-d} + c q^{d-2i} 
\qquad \qquad 0 \leq i \leq d,
\\
\label{eq:const2}
&&\theta^*_i = \alpha^* + b^* q^{2i-d} + c^* q^{d-2i} 
\qquad \qquad 0 \leq i \leq d,
\\
\label{eq:const3}
&&
q, 
\;\alpha,
\;b,
\;c,
\;\alpha^*,
\;b^*,
\;c^*
\; \in {\overline \F},
\\
\label{eq:const4}
&&
q \not=0, \quad q^2 \not=1,
\quad q^2 \not=-1.
\end{eqnarray}
We have a few comments about this solution.
For the moment define $\beta=q^2+q^{-2}$, and 
observe that $\beta+1$ is the common value of
(\ref{eq:betaplusone}).
We have $\beta-2=(q-q^{-1})^2$ and
 $\beta+2=(q+q^{-1})^2$. Therefore
 $\beta\not=2$, $\beta \not=-2$ in view of
(\ref{eq:const4}). Using
(\ref{eq:const1}),
(\ref{eq:const2}) we obtain
\begin{eqnarray*}
bc
&=&
\frac{
(\theta_0-\theta_1)^2
-\beta
(\theta_0-\theta_1)
(\theta_1-\theta_2)
+
(\theta_1-\theta_2)^2}{(\beta-2)^2(\beta+2)},
\\
b^*c^*
&=&
\frac{
(\theta^*_0-\theta^*_1)^2
-\beta
(\theta^*_0-\theta^*_1)
(\theta^*_1-\theta^*_2)
+
(\theta^*_1-\theta^*_2)^2}{
(\beta-2)^2(\beta+2)}
\end{eqnarray*}
provided $d\geq 2$.
We will focus on the case 
\begin{eqnarray}
bb^*cc^*\not=0.
\label{eq:const5}
\end{eqnarray}

\begin{definition}
\label{def:qr}
{\rm \cite[Definition~3.1]{IT:qRacah}}
\rm
Let $d$ denote a nonnegative integer and
let 
$(\lbrace \theta_i\rbrace_{i=0}^d;
\lbrace \theta^*_i\rbrace_{i=0}^d)$ 
denote a sequence of scalars taken from
$\F$. We call this sequence
{\it $q$-Racah}
whenever the following {\rm  (i), (ii)} hold:
\begin{enumerate}
\item[\rm (i)]
$\theta_i \neq \theta_j$, 
$\theta^*_i \neq \theta^*_j$ if $i \neq j$ $(0 \leq i,j \leq d)$;
\item[\rm (ii)]
there exist 
$q, \alpha,b,c,\alpha^*,b^*,c^*$ 
that satisfy
(\ref{eq:const1})--(\ref{eq:const5}).
\end{enumerate}
\end{definition}

\begin{definition}
\label{def:RAC}
\rm
For all integers $d\geq 0$ let
${\rm Rac}(d,\F)$ 
denote the set of all
 $q$-Racah
 sequences
 $(\lbrace \theta_i\rbrace_{i=0}^d;
\lbrace \theta^*_i\rbrace_{i=0}^d)$ 
of scalars taken from $\F$.
\end{definition}

\noindent
Observe that the set
${\rm Rac}(d,\F)$ 
from
Definition
\ref{def:RAC}
is contained in the set
${\rm Feas}(d,\F)$ 
from
Definition
\ref{def:feaset}.
In the next lemma
we characterize 
${\rm Rac}(d,\F)$ as a subset of
${\rm Feas}(d,\F)$.
To avoid trivialities we assume $d\geq 3$.

\begin{lemma}
\label{lem:qrrec}
Fix an integer $d\geq 3$ and
 a sequence $(\lbrace \theta_i\rbrace_{i=0}^d;
\lbrace \theta^*_i\rbrace_{i=0}^d)$ 
in
${\rm Feas}(d,\F)$.
Let
$\beta+1$ denote the common value of
(\ref{eq:betaplusone}).
Then the sequence  
is
in
${\rm Rac}(d,\F)$
 if and only if each of the following hold:
\begin{enumerate}
\item[\rm (i)] $\beta^2 \not=4$;
\item[\rm (ii)] 
$(\theta_0-\theta_1)^2
-\beta
(\theta_0-\theta_1)
(\theta_1-\theta_2)
+
(\theta_1-\theta_2)^2
\not=0$;
\item[\rm (iii)] 
$(\theta^*_0-\theta^*_1)^2
-\beta
(\theta^*_0-\theta^*_1)
(\theta^*_1-\theta^*_2)
+
(\theta^*_1-\theta^*_2)^2
\not=0$.
\end{enumerate}
\end{lemma}
\noindent {\it Proof:}
Use the comments above Definition
\ref{def:qr}.
\hfill $\Box$ \\

\begin{proposition}
\label{lem:bootstrap}
Assume $\F$ is infinite
and pick an integer $d\geq 3$.
Let $h$ denote a polynomial in
$2d+2$ mutually commuting indeterminates
that has all coefficients in $\F$.
Suppose that $h(p)=0$ for all
$p \in 
{\rm Rac}(d,\F)$.
Then
$h(p)=0$ for all
$p \in 
{\rm Feas}(d,\F)$.
\end{proposition}
\noindent {\it Proof:}
Let $\flat,
\lbrace t_i\rbrace_{i=0}^2,
\lbrace t^*_i\rbrace_{i=0}^2$
denote mutually commuting indeterminates.
Consider the $\F$-algebra 
$\F\lbrack 
\flat, t_0, t_1, t_2, t^*_0, t^*_1, t^*_2\rbrack$ 
consisting of the polynomials
in 
$\flat, t_0, t_1, t_2, t^*_0, t^*_1, t^*_2$
that have all coefficients in $\F$.
For $3 \leq i\leq d$ define $t_i, t^*_i \in 
\F\lbrack 
\flat, t_0, t_1, t_2, t^*_0, t^*_1, t^*_2\rbrack$ 
 by 
\begin{eqnarray*}
0 &=& t_i-(\flat+1)t_{i-1}+(\flat+1)t_{i-2}-t_{i-3},
\\
0 &=& t^*_i-(\flat+1)t^*_{i-1}+(\flat+1)t^*_{i-2}-t^*_{i-3}.
\end{eqnarray*}
Define a polynomial $f \in 
\F\lbrack 
\flat, t_0, t_1, t_2, t^*_0, t^*_1, t^*_2\rbrack$
to be the composition
\begin{eqnarray*}
f&=&
h(t_0,t_1,\ldots, t_d,t^*_0,t^*_1,\ldots,t^*_d).
\end{eqnarray*}
We mention one significance of $f$.
Given a sequence $p=(\lbrace \theta_i\rbrace_{i=0}^d;
\lbrace \theta^*_i\rbrace_{i=0}^d)$ in
${\rm Feas}(d,\F)$, 
let $\beta+1$ denote the common value of
(\ref{eq:betaplusone}). 
Define the sequence 
$s=(\beta,\theta_0,\theta_1,\theta_2,\theta^*_0,\theta^*_1,\theta^*_2)$.
Observe that 
$\theta_i = 
t_i(s) $ and
$\theta^*_i = 
t^*_i(s)$
for
$0 \leq i \leq d$.
Therefore
\begin{eqnarray}
\label{eq:fg}
f(s)=
h(p).
\end{eqnarray}
We show $f=0$.
Instead of working directly with $f$, it will be convenient
to work with the product
$\Psi=f\xi \xi^* \omega \omega^*(\flat^2-4)$, 
where 
\begin{eqnarray}
\label{eq:prod1}
\xi&=&\prod_{0 \leq i <j\leq d}
(t_i-t_j),
\\
\label{eq:prod2}
\xi^*&=&
\prod_{0 \leq i <j\leq d}
(t^*_i-t^*_j),
\\
\label{eq:prod3}
\omega &=&(t_0-t_1)^2
-\flat
(t_0-t_1)
(t_1-t_2)
+
(t_1-t_2)^2,
\\
\label{eq:prod4}
\omega^*&=&(t^*_0-t^*_1)^2
-\flat
(t^*_0-t^*_1)
(t^*_1-t^*_2)
+
(t^*_1-t^*_2)^2.
\end{eqnarray}
Each of 
$\xi, \xi^*$ is nonzero
by Lemma
\ref{lem:pexist} and since $\F$ is infinite.
Each of
$\omega,\omega^*,\flat^2-4 $ is nonzero by construction.
To show $f=0$ we will show that
$\Psi=0$ and invoke the fact that
$\F\lbrack 
\flat, t_0, t_1, t_2, t^*_0, t^*_1, t^*_2\rbrack$ is a domain
\cite[p.~129]{rotman}.
We now show that
$\Psi=0$.
Since  $\F$ is infinite it suffices 
to show that 
$\Psi(s)=0$
for all sequences 
$s=(\beta,\theta_0,\theta_1,\theta_2,\theta^*_0,\theta^*_1,\theta^*_2)$
of scalars
taken from $\F$
\cite[Proposition~6.89]{rotman}.
 Let $s$ be given.
For 
 $3 \leq i \leq d$ define
$\theta_i = t_i(s)$,
$\theta^*_i = t^*_i(s)$ 
and put 
 $p=(\lbrace \theta_i\rbrace_{i=0}^d;
\lbrace \theta^*_i\rbrace_{i=0}^d)$. 
We may assume
 $\lbrace \theta_i\rbrace_{i=0}^d$ are mutually distinct;
otherwise
$\xi(s)=0$
so
$\Psi(s)=0$.
We may assume
 $\lbrace \theta^*_i\rbrace_{i=0}^d$ are mutually distinct;
otherwise
$\xi^*(s)=0$
so
$\Psi(s)=0$.
By construction $p$ satisfies
condition (ii) of
Theorem \ref{conj:main},
with $\beta+1$ the  common value of 
(\ref{eq:betaplusone}).
Therefore $p$ is feasible by Definition
\ref{def:feas}.
For the moment assume that $p \in
{\rm Rac}(d,\F)$.
Then 
$f(s)=0$
by
(\ref{eq:fg}) and since
$h(p)=0$.
Therefore
$\Psi(s)=0$.
Next assume that
$ p \not\in {\rm Rac}(d,\F)$.
Then the product $\omega \omega^*(\flat^2-4)$ vanishes at
$s$
in view of Lemma
\ref{lem:qrrec}. 
The product 
$\omega \omega^*(\flat^2-4)$ 
is a factor of $\Psi$ so
$\Psi(s)=0$.
By the above comments 
$\Psi(s)=0$
for all sequences of scalars
$s=
(\beta,\theta_0,\theta_1,\theta_2,\theta^*_0,\theta^*_1,\theta^*_2)$
taken from $\F$. 
Therefore $\Psi=0$  so
$f=0$.
Now consider any sequence 
 $p=(\lbrace \theta_i\rbrace_{i=0}^d;
\lbrace \theta^*_i\rbrace_{i=0}^d)$
in
${\rm Feas}(d,\F)$.
 Then $h(p)=0$ by
(\ref{eq:fg}) and since $f=0$.
\hfill $\Box$ \\

\section{The algebra $\check T$}

\noindent In order to prove
Theorem
\ref{conj:mainp}  
we will need some detailed results about the
algebra $T$
from Definition
\ref{def:twitha}. In order to obtain
these results 
it is helpful to 
first consider 
the following 
algebra $\check T$.

\begin{definition}
\rm
\label{def:free}
Fix an integer $d\geq 0$.
Let ${\check T}={\check T}(d,\F)$ denote the $\F$-algebra defined by
generators 
$\lbrace \epsilon_i\rbrace_{i=0}^d$,
$\lbrace \epsilon^*_i\rbrace_{i=0}^d$
and  relations 
\begin{eqnarray}
\label{eq:epsrel}
&&\epsilon_i
\epsilon_j = \delta_{i,j} \epsilon_i,
\qquad \qquad
\epsilon^*_i\epsilon^*_j = \delta_{i,j} \epsilon^*_i
\qquad \qquad 0 \leq i,j\leq d.
\end{eqnarray}
\end{definition}

\begin{definition}
\label{def:igen}
\rm
Referring to Definition 
\ref{def:free},
we call 
 $\lbrace \epsilon_i\rbrace_{i=0}^d$ and
 $\lbrace \epsilon^*_i\rbrace_{i=0}^d$ the {\it idempotent 
generators} for $\check T$.
We say that the  
 $\lbrace \epsilon^*_i\rbrace_{i=0}^d$ are {\it starred}
and the 
 $\lbrace \epsilon_i\rbrace_{i=0}^d$ are {\it nonstarred}.
\end{definition}

\begin{definition}
\label{def:wordfree}
\rm
A pair of idempotent generators for $\check T$
is called {\it alternating} whenever one
of them is starred and the other is nonstarred.
For an integer $n\geq 0$,
by a {\it word of length $n$} in $\check T$
we mean a product $g_1g_2\cdots g_n$ such that
$\lbrace g_i\rbrace_{i=1}^n$
are idempotent generators for $\check T$
and $g_{i-1},g_i$ are alternating for $2 \leq i\leq n$.
We interpret the word of length $0$ to be the 
identity of $\check T$. We call this word
{\it trivial}.
\end{definition}

\begin{proposition}
\label{prop:basisfree}
The $\F$-vector space $\check T$ has a basis
consisting of its words.
\end{proposition}
\noindent {\it Proof:}
Let $S$ denote the set of words in $\check T$.
By construction $S$ spans $\check T$.
We show that $S$ is linearly independent.
To this end we introduce 
some indeterminates
$\lbrace f_i\rbrace_{i=0}^d$;
$\lbrace f^*_i\rbrace_{i=0}^d$ called
{\it formal idempotents}. We call the
$\lbrace f^*_i\rbrace_{i=0}^d$ {\it starred}
and the 
$\lbrace f_i\rbrace_{i=0}^d$ {\it nonstarred}.
A pair of formal idempotents is said to be
{\it alternating} whenever one of them is starred
and the other is nonstarred. For an integer
$n\geq 0$, by a {\it formal word of length $n$}
we mean a sequence $(y_1, y_2, \ldots, y_n)$ such that
$\lbrace y_i\rbrace_{i=1}^n$ are formal idempotents
and $y_{i-1},y_i$ are alternating  for $2 \leq i \leq n$.
The formal word of length 0 is called {\it trivial}
and denoted by $1$. Let $\mathcal S$
denote the set of all formal words. Let $V$
denote the vector space over $\F$ consisting of the
$\F$-linear combinations of $\mathcal S$ that have finitely
many nonzero coefficients. The set $\mathcal S$
is a basis for $V$. For $0 \leq i \leq d$ we define
linear transformations
$F_i: V\to V$ and
$F^*_i : V\to V$.
To do this we give the action of
$F_i$ and
$F^*_i$ on $\mathcal S$. We define
$F_i.1 = f_i$
and 
$F^*_i.1 = f^*_i$.
Pick a nontrivial formal word
$y=(y_1,y_2,\ldots, y_n)$. For the moment 
assume that $y_1$ is starred.
We define
$F_i.y=
(f_i,y_1,y_2,\ldots, y_n)$.
Also 
$F^*_i.y=y$ if $y_1=f^*_i$ and
$F^*_i.y=0$ if $y_1\not=f^*_i$.
Next assume 
that $y_1$ is nonstarred.
We define
$F_i.y=y$ if $y_1=f_i$ and
$F_i.y=0$ if $y_1\not=f_i$.
Also
$F^*_i.y=
(f^*_i,y_1,y_2,\ldots, y_n)$.
The linear transformations
$F_i:V\to V$ and $F^*_i:V \to V$ are now defined.
 By construction
$F_iF_j = \delta_{i,j}F_i$ and
$F^*_iF^*_j = \delta_{i,j}F^*_i$ for
$0 \leq i,j\leq d$.
Therefore $V$ has a $\check T$-module structure such that
$\epsilon_i$ (resp. $\epsilon^*_i$)
acts on $V$ as $F_i$ (resp. $F^*_i$)
for $0 \leq i \leq d$.
Consider the linear transformation $\gamma:{\check T}\to V$
that sends
$z \mapsto z.1$ for all $z \in {\check T}$.
For each word $g_1g_2\cdots g_n$ in $\check T$ we find
$\gamma(g_1g_2\cdots g_n)=(g'_1,g'_2, \ldots, g'_n)$,
where $\epsilon'_i=f_i$ and
$\epsilon^{*\prime}_i=f^*_i$ for $0 \leq i \leq d$.
Therefore the restriction of $\gamma$ to
$S$ gives a bijection 
$S\mapsto {\mathcal S}$. The set
${\mathcal S}$ is linearly independent
and $\gamma$ is linear so $S$ is linearly independent.
We have shown that $S$ is a basis for $\check T$.
\hfill $\Box$ \\

\noindent 
Let $u,v$ denote words in $\check T$.
Then their product $uv$ is either 0 or a word in $\check T$.

\section{The algebras $ D$
and $D^*$}

\noindent Throughout this section
we fix
an integer $d\geq 0$  and consider the algebra
${\check T}={\check T}(d,\F)$
from Definition
\ref{def:free}.

\begin{definition}
\label{def:d} 
\rm
Let $D$ (resp. $D^*$) denote the
subspace of
${\check T}$ with a basis
$\lbrace \epsilon_i\rbrace_{i=0}^d$
(resp. 
$\lbrace \epsilon^*_i\rbrace_{i=0}^d$).
\end{definition}

\noindent  We mention some  notation.
 For subsets $Y, Z$ of $\check T$
let $YZ$ denote the subspace of $\check T$ spanned
by $\lbrace yz \, |\, y \in Y, \; z \in Z\rbrace$. 

\begin{lemma}
\label{lem:ddst}
In the $\F$-vector space $\check T$ the
following sum is direct:
\begin{eqnarray*}
{\check T} = \F 1 \;+\; D \; + \;D^* \;+ \; DD^* \; + \; D^*D
\;+ \;DD^*D \; + \; D^*DD^* \; + \; \cdots 
\end{eqnarray*}
Moreover the $\F$-algebra $\check T$ is generated by $D,D^*$.
\end{lemma}
\noindent {\it Proof:} The first assertion is immediate from Proposition
\ref{prop:basisfree}. The last assertion is clear.
\hfill $\Box$ \\

\begin{lemma}
\label{lem:ddsident}
The space $D$ (resp. $D^*$) 
is an $\F$-algebra
with multiplicative
identity
$\sum_{i=0}^d \epsilon_i $
(resp. $\sum_{i=0}^d \epsilon^*_i $).
\end{lemma}
\noindent {\it Proof:} 
Use 
(\ref{eq:epsrel}).
\hfill $\Box$ \\


\noindent 
We emphasize that $D$ and $D^*$ are not subalgebras of
$\check T$,  since their multiplicative identities  do not
equal the multiplicative identity 1 of $\check T$.
However we do have the following.

\begin{lemma}
\label{lem:subalgd}
The spaces  
$D+\F 1$ and
$D^*+\F 1 $ 
are $\F$-subalgebras of $\check T$.
\end{lemma}
\noindent {\it Proof:}
These subspaces
are  closed under multiplication
and contain the identity $1$
of $\check T$.
\hfill $\Box$ \\

\section{The homogeneous components of $\check T$}

\noindent Throughout this section
we fix
an integer $d\geq 0$  and consider the algebra
${\check T}={\check T}(d,\F)$
from Definition
\ref{def:free}.

\begin{definition}
\label{def:beginend}
\rm
Let $w=g_1g_2\cdots g_n$ denote a nontrivial word in $\check T$.
We say that $w$ {\it begins} with $g_1$ and
{\it ends} with $g_n$. We write
\begin{eqnarray*}
g_1 = {\rm begin}(w),
\qquad \qquad 
g_n = {\rm end}(w).
\end{eqnarray*}
\end{definition}

\begin{example}
\rm
Assume $d=2$.
In the table below we display some nontrivial words
$w$ in $\check T$. For each word $w$
we give ${\rm begin}(w)$ and
${\rm end}(w)$.

\begin{center}
\begin{tabular}{c|c c}
$w$ &  ${\rm begin}(w)$ & ${\rm end}(w)$  \\
\hline \hline
$\epsilon_1$  & $\epsilon_1$ & $\epsilon_1$ 
\\
$\epsilon_1 \epsilon^*_2$ & $\epsilon_1$ & $\epsilon^*_2$ 
\\
$\epsilon^*_1 \epsilon_0 \epsilon^*_2$  & $\epsilon^*_1$ & $\epsilon^*_2$ 
    \end{tabular}
    \end{center}
\end{example}

\begin{definition}
\label{def:sim}
\rm
We define a binary
relation $\sim$ on the set of words in $\check T$.
With respect to $\sim$ 
the trivial word in $\check T$ is related to itself and no other
word in $\check T$.
For nontrivial words $u, v$ in $\check T$
we define
$u\sim v$ whenever
each of the following holds:
\begin{eqnarray*}
{\rm length}(u)=
{\rm length}(v),
\qquad \qquad 
{\rm begin}(u)=
{\rm begin}(v),
\qquad \qquad 
{\rm end}(u)=
{\rm end}(v).
\end{eqnarray*}
Observe that $\sim$ is an equivalence relation.
\end{definition}


\begin{definition}
\label{def:eqclasses}
\rm
Let $\Lambda $ denote
the set of equivalence classes for the relation
$\sim$ in
 Definition
\ref{def:sim}.
An element of $\Lambda$ is called a {\it type}.
For $\lambda \in \Lambda$ the words in $\lambda$
are said to have {\it type $\lambda$}.
\end{definition}

\begin{definition}
\label{def:typelength}
\rm
For $\lambda \in \Lambda$ let 
${\rm length}(\lambda)$ denote the common length of
each word of type $\lambda$.
\end{definition}

\begin{definition}
\label{def:trivialtype}
\rm 
There exists a 
unique type
in $\Lambda$ that has length $0$. This type consists
of the trivial word $1$ and nothing else.
We call
this type {\it trivial}.
\end{definition}

\begin{definition}
\rm
For all nontrivial $\lambda \in \Lambda$,
\begin{enumerate}
\item[\rm (i)]
let ${\rm begin}(\lambda)$ denote the common 
beginning of each word of type $\lambda$;
\item[\rm (ii)] 
let ${\rm end}(\lambda)$ denote the common 
ending of each word of type $\lambda$.
\end{enumerate}
\end{definition}

\begin{definition}
\label{def:tlambda}
\rm
For $\lambda \in \Lambda$ let ${\check T}_{\lambda}$
denote the subspace of $\check T$ 
with a basis
consisting of the words of type $\lambda$.
\end{definition}

\begin{proposition}
\label{thm:dec}
The $\F$-vector space 
${\check T}$
decomposes as
\begin{eqnarray}
{\check T} = \sum_{\lambda \in \Lambda}
{\check T}_{\lambda}
\qquad \qquad (\mbox{\rm direct sum}).
\label{eq:tcheckds}
\end{eqnarray}
\end{proposition}
\noindent {\it Proof:}
Immediate from
Proposition
\ref{prop:basisfree}
and 
Definition \ref{def:tlambda}.
\hfill $\Box$ \\


\begin{definition}
\rm
For $\lambda \in \Lambda$ we call
${\check T}_{\lambda}$ the {\it $\lambda$-homogeneous component}
of $\check T$.
Elements of 
${\check T}_{\lambda}$
are said to be {\it $\lambda$-homogeneous}.
An element  of 
${\check T}$
is called {\it homogeneous} whenever 
it is $\lambda$-homogeneous for some $\lambda \in \Lambda$.
\end{definition}


\section{The zigzag words in $\check T$}

\noindent Throughout this section
we fix
an integer $d\geq 0$  and consider the algebra
${\check T}={\check T}(d,\F)$
from Definition
\ref{def:free}.
We have been discussing the
words in 
 $\check T$.
We now focus our attention on a
special kind of word
 said to be zigzag.

\begin{definition}
\rm 
Given an ordered pair of integers
$i,j$ and an integer $m$
we say that
$m$ is {\it between}
$i,j$ whenever $i\geq m>j$
or 
$i\leq m<j$.
\end{definition}

\begin{definition}
\label{def:indexf}
\rm
For an idempotent generator $\epsilon_i$ or
$\epsilon^*_i$ of $\check T$, we call $i$ the {\it index} of the 
generator. For an idempotent generator $g$
of $\check T$
let $\overline g$ denote the index of $g$.
\end{definition}

\begin{definition}
\label{def:zzf}
\rm
A word $g_1g_2\cdots g_n$ in $\check T$ is said
to be {\it zigzag} whenever both
\begin{enumerate}
\item[\rm (i)] ${\overline g}_i$ is not between 
${\overline g}_{i-1}$,
${\overline g}_{i+1}$ for $2 \leq i \leq n-1$;
\item[\rm (ii)] at least one of
${\overline g}_{i-1}$,
${\overline g}_i$ is not between 
${\overline g}_{i-2}$,
${\overline g}_{i+1}$ for $3 \leq i \leq n-1$.
\end{enumerate}
\end{definition}

\noindent  We now describe 
the zigzag words in $\check T$.
We will use the following notion.
Two integers $m,m'$ are said to have {\it opposite
sign} whenever $m m'\leq 0$.

\begin{proposition}  
{\rm \cite[Theorem~7.7]{NT:onshape}}
\label{thm:zzz}
Let $g_1g_2\cdots g_n$ denote a word in $\check T$.
Then this word is zigzag if and only if both
\begin{itemize}
\item[\rm (i)]
${\overline g}_{i-1}-{\overline g}_{i}$
and ${\overline g}_{i}-{\overline g}_{i+1}$ have opposite sign
for $2 \leq i \leq n-1$;
\item[\rm (ii)]
for $2 \leq i \leq n-1$, if $\;|{\overline g}_{i-1}-{\overline g}_{i}|
<|{\overline g}_{i}-{\overline g}_{i+1}|$
then
\[
 0 < |{\overline g}_1-{\overline g}_2| < |{\overline g}_2-{\overline g}_3|
 < \cdots < |{\overline g}_{i}-{\overline g}_{i+1}|.
 \]
 \end{itemize}
 \end{proposition}

\begin{definition}    \label{def:pq}
\rm
A word $g_1g_2\cdots g_n$ in $\check T$ 
is said to be {\it constant} whenever the index
${\overline g}_i$ is independent of $i$ for $1 \leq i \leq n$.
Note that the trivial word is constant, and each constant word is zigzag.
\end{definition}

\begin{proposition} 
{\rm \cite[Theorem~7.9]{NT:onshape}}
\label{thm:p} 
Let $g_1g_2\cdots g_n$ denote a nonconstant zigzag word in
$\check T$.
Then there exists a unique integer $\kappa $ $(2 \leq \kappa \leq n)$ 
such that both
\begin{itemize}
\item[\rm (i)]
 $0 < |{\overline g}_1-{\overline g}_2| < \cdots < 
 |{\overline g}_{\kappa-1}-{\overline g}_{\kappa}|$;
 \item[\rm (ii)]
  $|{\overline g}_{\kappa-1}-{\overline g}_{\kappa}| \geq 
  |{\overline g}_{\kappa}-{\overline g}_{\kappa+1}| \geq \cdots \geq 
  |{\overline g}_{n-1}-{\overline g}_n|$.
  \end{itemize}
  \end{proposition}

\begin{definition}
\label{def:zlambda}
\rm
For $\lambda \in \Lambda$ let $Z_{\lambda}$
denote the subspace of ${\check T}$
with a basis consisting of the zigzag words of type $\lambda$.
Note that $Z_{\lambda} \subseteq {\check T}_{\lambda}$.
\end{definition}

\section{The algebra $\epsilon^*_0 {\check T}\epsilon^*_0$}

\noindent Throughout this section
we fix
an integer $d\geq 0$  and consider the algebra
${\check T}={\check T}(d,\F)$
from Definition
\ref{def:free}.
 Observe that $\epsilon^*_0{\check T}\epsilon^*_0$
is an $\F$-algebra with multiplicative identity 
$\epsilon^*_0$.

\begin{lemma}
\label{lem:tzerobasis}
The $\F$-vector space 
$\epsilon^*_0{\check T}\epsilon^*_0$
has a basis consisting of the nontrivial words in $\check T$ that
begin and end with $\epsilon^*_0$.
\end{lemma}
\noindent {\it Proof:} 
Let $U$ denote the
subspace of $\check T$ with a basis consisting of
the nontrivial words in $\check T$
that begin and end with $\epsilon^*_0$.
We show that
$\epsilon^*_0{\check T}\epsilon^*_0=U$.
We first show that
$\epsilon^*_0{\check T}\epsilon^*_0 \subseteq U$.
Recall that $\check T$ is spanned by its words.
For all words $w$ in
${\check T}$ the product
$\epsilon^*_0 w\epsilon^*_0$
is either zero, or a nontrivial word in $\check T$
that begins and ends with 
 $\epsilon^*_0$.
In either case
$\epsilon^*_0 w\epsilon^*_0 \in U$, and therefore
$\epsilon^*_0{\check T}\epsilon^*_0 \subseteq U$.
Next we show that
$U \subseteq \epsilon^*_0{\check T}\epsilon^*_0$.
Let $w$ denote a nontrivial word 
in $\check T$ that begins and ends with $\epsilon^*_0$.
We have $w=\epsilon^*_0 w \epsilon^*_0$ since
$\epsilon^{*2}_0 =
 \epsilon^*_0$, so
$ w \in \epsilon^*_0{\check T}\epsilon^*_0$.
Therefore $U \subseteq \epsilon^*_0{\check T}\epsilon^*_0$.
We have shown that
 $\epsilon^*_0{\check T}\epsilon^*_0=U$ and the result follows.
\hfill $\Box$ \\

\begin{definition}
\label{def:lambdazero}
\rm
Let $\Lambda_0$ denote the set of types
in $\Lambda$ that begin and end with $\epsilon^*_0$.
\end{definition}

\noindent Our next goal is to describe
$\Lambda_0$.

\begin{definition}
\label{def:starlength}
\rm
Let $g_1g_2\cdots g_n$ denote a
word in $\check T$. 
By the {\it star-length} (resp. {\it nonstar-length}) of
this word we mean the number of
terms in the sequence
$(g_1, g_2,\ldots, g_n)$ that are
starred (resp. nonstarred).
Note that the star-length plus the nonstar-length is
equal to the length $n$.
For $\lambda \in \Lambda$,
by the 
 {\it star-length} (resp. {\it nonstar-length}) of
$\lambda$
we mean the common 
 star-length (resp. nonstar-length) of
each word of type $\lambda$.
\end{definition}

\begin{definition}
\label{def:lambdazero2}
\rm
For an integer  $n\geq 0$ let 
$\lbrack n\rbrack$ denote the
unique type in 
$\Lambda_0$ that has nonstar-length $n$. 
Observe that
$\lbrack n\rbrack$ has star-length $n+1$
and length $2n+1$.
\end{definition}

\noindent The next two lemmas follow immediately from
Definition
\ref{def:lambdazero}
and Definition
\ref{def:lambdazero2}.

\begin{lemma}
\label{lem:biject}
The map 
$n \mapsto \lbrack n \rbrack$ gives a bijection
from the set of nonnegative integers to the set
$\Lambda_0$.
\end{lemma}

\begin{lemma}
\label{lem:nmult}
Let $m$ and $n$ denote nonnegative integers.
Let $u$ and $v$ denote words in $\check T$
of type $\lbrack m\rbrack $ and
$\lbrack n\rbrack $ respectively.
Then $uv$ is a word in
 $\check T$
of type $\lbrack m+n\rbrack $.
\end{lemma}

\begin{proposition}
\label{lem:ete1}
The $\F$-vector space 
$\epsilon^*_0{\check T}\epsilon^*_0$ decomposes as
\begin{eqnarray}
\label{eq:tzerodec}
\epsilon^*_0{\check T}\epsilon^*_0
=
\sum_{n=0}^\infty
{\check T}_{\lbrack n \rbrack}
\qquad \qquad  (\mbox{\rm direct sum}).
\end{eqnarray}
Moreover 
${\check T}_{\lbrack m \rbrack}
\cdot
{\check T}_{\lbrack n \rbrack}
\subseteq 
{\check T}_{\lbrack m+n \rbrack}
$
for all integers $m,n\geq 0$.
\end{proposition}
\noindent {\it Proof:}
By Lemma
\ref{lem:tzerobasis}
and Definition
\ref{def:lambdazero}
we have
$\epsilon^*_0{\check T}\epsilon^*_0
=\sum_{\lambda \in \Lambda_0}
{\check T}_\lambda$
(direct sum).
Combining
this with
Lemma
\ref{lem:biject}
we obtain
(\ref{eq:tzerodec}).
The last assertion follows from
Lemma
\ref{lem:nmult}.
\hfill $\Box$ \\


\noindent We turn our attention to
the zigzag words in $\check T$ that
 begin and end with $\epsilon^*_0$.

\begin{proposition}
\label{prop:zzn}
Pick an integer $n\geq 0$ and a word
$g_1g_2\cdots g_{2n+1}$ in $\check T$ of type 
$\lbrack n \rbrack$.
This word is zigzag if and only if
both 
\begin{enumerate}
\item[\rm (i)]
${\overline g}_i=0$ for all odd $i$ $(1 \leq i \leq 2n+1)$;
\item[\rm (ii)]
${\overline g}_i\geq 
{\overline g}_{i+2}
$
 for all even $i$ $(2 \leq i \leq 2n-2)$.
\end{enumerate}
\end{proposition}
\noindent {\it Proof:}
The type 
$\lbrack n \rbrack$ begins and ends with
$\epsilon^*_0$, so
${\overline g}_1=0$ and
${\overline g}_{2n+1}=0$.
The result follows from 
this and
Propositions
\ref{thm:zzz},
\ref{thm:p}. 
\hfill $\Box$ \\

\section{The elements $a$, $a^*$}

\medskip
\noindent 
Recall that the algebra $T$ from
Definition
\ref{def:twitha} is defined using
relations
(\ref{eq:eiejcop})--(\ref{eq:eiaskej}).
So far we have investigated relation
(\ref{eq:eiejcop}). We now
prepare to bring in relations
(\ref{eq:sumthieicop})--(\ref{eq:eiaskej}).

\medskip
\noindent Throughout this section
we fix an integer  $d \geq 0$ 
and 
 a sequence 
$p=(\lbrace \theta_i\rbrace_{i=0}^d;
\lbrace \theta^*_i\rbrace_{i=0}^d)$ 
in 
${\rm Feas}(d,\F)$. 
Recall the
algebra ${\check T}={\check T}(d,\F)$
from Definition
\ref{def:free}. 

\begin{definition}
\label{def:adef}
\rm Define $a=a(p)$ and
$a^* =a^*(p)$ in $\check T$ by
\begin{eqnarray}
\label{eq:adef}
a = \sum_{i=0}^d \theta_i \epsilon_i,
\qquad \qquad
a^* = \sum_{i=0}^d \theta^*_i \epsilon^*_i.
\end{eqnarray}
Observe that $a \in D$ and $a^* \in D^*$, where $D$, $D^*$ are
from Definition
\ref{def:d}.
\end{definition}

\begin{lemma}
For $0 \leq i \leq d$,
\begin{eqnarray}
&&\quad a\epsilon_i=\epsilon_i a = \theta_i \epsilon_i,
\qquad \qquad 
a^*\epsilon^*_i =\epsilon^*_i a^* =  \theta^*_i \epsilon^*_i.
\label{eq:aprod}
\end{eqnarray}
\end{lemma}
\noindent {\it Proof:}
Use
(\ref{eq:epsrel})
and
(\ref{eq:adef}).
\hfill $\Box$ \\

\begin{note}
\label{ameaning}
\rm  
We will be considering powers of
the elements $a$, $a^*$ from
Definition
\ref{def:adef}.
We wish to clarify the meaning
of $a^0$ and $a^{*0}$.
We always interpret
\begin{eqnarray}
a^0 = \sum_{i=0}^d \epsilon_i, \qquad \qquad
a^{*0}=\sum_{i=0}^d \epsilon^*_i.
\label{eq:meaning}
\end{eqnarray}
This  is justified by
Lemma
\ref{lem:ddsident}.
We mention some related notational conventions.
Consider the $\F$-algebra homomorphism
$\F\lbrack x \rbrack \to D$ 
that sends
$x \mapsto a$.
By definition this homomorphism sends
the identity $1$ of 
 $\F\lbrack x \rbrack$ to the identity
 $\sum_{i=0}^d \epsilon_i $ of
 $D$.
For $f \in \F\lbrack x \rbrack$ the image
of $f$ under this homomorphism will be denoted $f(a)$.
Writing $f=\sum_{i=0}^n c_i x^i$ we have
$f(a) = \sum_{i=0}^n c_i a^i$,
with the $i=0$ summand interpreted using
the equation on the left in
(\ref{eq:meaning}). A similar comment applies to
$a^*$.
\end{note}
\noindent The above notational conventions are illustrated
in the following lemma.

\begin{lemma}
\label{lem:apower} For $f\in \F\lbrack x \rbrack$,
\begin{eqnarray*}
f(a)= \sum_{i=0}^d f(\theta_i)\epsilon_i,
\qquad\qquad
f(a^*)= \sum_{i=0}^d f(\theta^*_i)\epsilon^*_i.
\end{eqnarray*}
In particular
for an integer $k\geq 0$,
\begin{eqnarray}
a^k = \sum_{i=0}^d \theta^k_i \epsilon_i,
\qquad \qquad 
a^{*k} = \sum_{i=0}^d \theta^{*k}_i \epsilon^*_i.
\label{eq:apower}
\end{eqnarray}
\end{lemma}
\noindent {\it Proof:}
Use
(\ref{eq:epsrel})
and (\ref{eq:adef}).
\hfill $\Box$ \\

\section{The algebra ${\tilde T}$}

\noindent In our study of the algebra $T$
we now
bring in relations
(\ref{eq:sumthieicop})--(\ref{eq:eiaskej}).
We do this in a compact way.

\begin{definition}
\label{def:tt}
\rm
Fix an integer $d \geq 0$  and a sequence
 $p=(\lbrace \theta_i\rbrace_{i=0}^d;
\lbrace \theta^*_i\rbrace_{i=0}^d)$ 
in
${\rm Feas}(d,\F)$.
Let ${\tilde T}={\tilde T}(p,\F)$ denote the $\F$-algebra with generators
$\lbrace \varepsilon_i\rbrace_{i=0}^d$,
$\lbrace \varepsilon^*_i\rbrace_{i=0}^d$ and relations
\begin{eqnarray}
&&\varepsilon_i\varepsilon_j = \delta_{i,j} \varepsilon_i,
\qquad \qquad
\varepsilon^*_i\varepsilon^*_j = \delta_{i,j} \varepsilon^*_i
\qquad \qquad 0 \leq i,j\leq d,
\label{eq:r1}
\\
&&0=\sum_{\ell=0}^d \theta^k_\ell \varepsilon^*_i \varepsilon_\ell
\varepsilon^*_j,
\qquad
0=\sum_{\ell=0}^d \theta^{*k}_\ell \varepsilon_i
\varepsilon^*_\ell \varepsilon_j
 \qquad 0 \leq i,j\leq d,
\qquad 0 \leq k < |i-j|.
\qquad 
\label{eq:r2}
\end{eqnarray}
\end{definition}


\noindent Many of the concepts that apply to
$\check T$ also apply to $\tilde T$.
We emphasize a few such concepts in the following
definitions.

\begin{definition}
\label{def:ig}
\rm
Referring to Definition \ref{def:tt},
we call 
 $\lbrace \varepsilon_i\rbrace_{i=0}^d$ and
 $\lbrace \varepsilon^*_i\rbrace_{i=0}^d$ the {\it idempotent 
generators} for $\tilde T$.
We say that the  
 $\lbrace \varepsilon^*_i\rbrace_{i=0}^d$ are {\it starred}
and the 
 $\lbrace \varepsilon_i\rbrace_{i=0}^d$ are {\it nonstarred}.
A pair of idempotent generators for $\tilde T$
will be called {\it alternating} whenever one
of them is starred and the other is nonstarred.
\end{definition}

\begin{definition}
\label{def:word}
\rm
For an integer $n\geq 0$,
by a {\it word of length $n$} in $\tilde T$
we mean a product $g_1g_2\cdots g_n$ such that
$\lbrace g_i\rbrace_{i=1}^n$
are idempotent generators for $\tilde T$
and $g_{i-1},g_i$ are alternating for $2 \leq i\leq n$.
We interpret the word of length $0$ to be the 
identity of $\tilde T$. We call this word
{\it trivial}. Let $g_1g_2\cdots g_n$ denote a
nontrivial word in $\tilde T$. We say that this
word {\it begins} with $g_1$ and
{\it ends} with $g_n$.
\end{definition}

\noindent Referring to Definition
\ref{def:word}, observe that $\tilde T$ is
spanned by its words.

\medskip
\noindent From the construction we have 
canonical $\F$-algebra
homomorphisms 
${\check T} \to {\tilde T}
  \to  T$.
We will investigate these homomorphisms in
the following sections.

\section{The homomorphism $ \varphi: {\check T}\to {\tilde T}$}

\noindent  From now until
the end of Lemma
 \ref{lem:desirablemeaning}
the following notation
will be in effect. Fix an integer $d\geq 0$ and
let the algebra
${\check T}={\check T}(d,\F)$
be as in
Definition
\ref{def:free}.
 Fix a sequence  $p=(\lbrace \theta_i\rbrace_{i=0}^d;
\lbrace \theta^*_i\rbrace_{i=0}^d)$ 
in
${\rm Feas}(d,\F)$
 and let the algebra
${\tilde T}={\tilde T}(p,\F)$ be as
in
Definition
\ref{def:tt}.
We now consider how
${\check T}$ and ${\tilde T}$ are related.

\begin{definition}
\label{def:idealr}
\rm
Let $R=R(p)$ denote the two-sided ideal of $\check T$ generated
by the elements
\begin{eqnarray}
\label{eq:relvect}
\epsilon^*_i a^k \epsilon^*_j,
\qquad 
\epsilon_i a^{*k} \epsilon_j,
\qquad \qquad
0\leq i,j\leq d,\qquad 0 \leq k < |i-j|,
\end{eqnarray}
where $a=a(p)$ and $a^*=a^*(p)$ are from
Definition
\ref{def:adef}.
\end{definition}

\begin{lemma}
\label{lem:cancheck}
There exists a surjective $\F$-algebra homomorphism
$\varphi: {\check T}\to {\tilde T}$ that sends
$\epsilon_i \mapsto \varepsilon_i$ and
$\epsilon^*_i \mapsto \varepsilon^*_i$ for
$0 \leq i \leq d$.
The kernel of $\varphi$ coincides with the ideal $R$.
\end{lemma}
\noindent {\it Proof:}
Compare the defining relations for $\check T$ and
$\tilde T$.
\hfill $\Box$ \\

\noindent Our next goal is to display a spanning
set for $R$.
To this end we introduce a type of
element in $\check T$ called a relator.

\begin{definition}
\label{def:relatorpre}
\rm
Let $C$ (resp. $C^*$) denote the set of three-tuples
$(u,v,k)$ such that:
\begin{itemize}
\item[\rm (i)]
each of $u$, $v$ is a nontrivial word in $\check T$;
\item[\rm (ii)]
${\rm end}(u)$ and
${\rm begin}(v)$ are both nonstarred (resp. both starred);
\item[\rm (iii)] $k$ is an integer such that
$0 \leq k < |{\overline {{\rm end}(u)}}-{\overline {{\rm begin}(v)}}|$.
\end{itemize}
Observe that $C\cap C^* =\emptyset$.
\end{definition}

\begin{definition}
\label{def:relator}
\rm
With reference to
Definition \ref{def:relatorpre},
to each element in $C\cup C^*$ 
we associate an element of $\check T$
called its {\it relator}.
For $(u,v,k) \in C$ the corresponding
relator is
$ua^{*k} v$, where
$a^*=a^*(p)$ is from
Definition
\ref{def:adef}.
For $(u,v,k) \in C^*$
the corresponding relator is
$ua^{k} v$, where
$a=a(p)$ is from
Definition
\ref{def:adef}.
\end{definition}

\begin{lemma}
\label{lem:relatorspan}
The $\F$-vector space $R$ is spanned by
the relators in $\check T$.
\end{lemma}
\noindent {\it Proof:}
By Definition
\ref{def:idealr} and since $\check T$
is spanned by its words.
\hfill $\Box$ \\

\begin{lemma}
\label{lem:relatorspan2}
With reference to
Definition \ref{def:relatorpre},
for $(u,v,k) \in C\cup C^*$ the
corresponding relator is $\lambda$-homogeneous,
where
\begin{eqnarray*}
&&{\rm begin}(\lambda) = 
{\rm begin}(u),
\qquad \qquad
{\rm end}(\lambda) = 
{\rm end}(v),
\\
&&
\qquad {\rm length}(\lambda) = 
{\rm length}(u)+
{\rm length}(v)+1.
\end{eqnarray*}
\end{lemma}
\noindent {\it Proof:}
Use
(\ref{eq:apower}) and
Definition
\ref{def:relator}.
\hfill $\Box$ \\

\begin{definition}
\label{def:relatorsub}
\rm
For $\lambda \in \Lambda$ 
let 
$R_\lambda=
R_\lambda(p)
$ denote the subspace of
$\check T$ spanned by the 
$\lambda$-homogeneous relators.
Observe that 
$R_\lambda \subseteq 
{\check T}
_\lambda$.
\end{definition}

\begin{lemma}
\label{def:relvector2}
The $\F$-vector space $R$ 
decomposes as
\begin{eqnarray*}
R = \sum_{\lambda \in \Lambda} R_\lambda
\qquad \qquad    (\mbox{\rm direct sum}).
\end{eqnarray*}
\end{lemma}
\noindent {\it Proof:}
By Lemmas
\ref{lem:relatorspan},
\ref{lem:relatorspan2}
and Definition
 \ref{def:relatorsub}
we obtain 
$R = \sum_{\lambda \in \Lambda} R_\lambda $.
The sum
$\sum_{\lambda \in \Lambda} R_\lambda$
is direct by Proposition
\ref{thm:dec}
and since $R_\lambda \subseteq {\check T}_\lambda$
for all $\lambda \in \Lambda$.
\hfill $\Box$ \\

\begin{corollary}
\label{cor:rint}
For $\lambda \in \Lambda $ 
we have
$R_\lambda = R\cap {\check T}_\lambda$.
\end{corollary}
\noindent {\it Proof:}
Observe that
$R_\lambda \subseteq  R$
by Lemma
\ref{def:relvector2}
and
$R_\lambda \subseteq  {\check T}_\lambda$
by
Definition
\ref{def:relatorsub},
so 
$R_\lambda \subseteq R\cap {\check T}_\lambda$.
To obtain the reverse inclusion,
we pick any $v \in 
R\cap {\check T}_\lambda$ and show
$v \in 
R_\lambda $.
By Lemma
\ref{def:relvector2}
there exists 
$r \in R_\lambda$ such that
$v -r \in
\sum_{\chi \in \Lambda \backslash \lambda } R_\chi$.
We have
$v -r \in {\check T}_\lambda$
by construction and
the last sentence in Definition
\ref{def:relatorsub}.
Similarly $v-r
\in 
\sum_{\chi \in \Lambda \backslash \lambda } {\check T}_\chi$.
So 
$v -r$ is contained in the intersection of 
${\check T}_\lambda $ and
$\sum_{\chi \in \Lambda \backslash \lambda }{\check T}_\chi$.
Now 
$v =r$ in view of
Proposition
\ref{thm:dec}, so
$v \in  R_\lambda$.
We have shown $R_\lambda \supseteq R\cap {\check T}_\lambda$
and the result follows.
\hfill $\Box$ \\

\begin{definition}
\label{def:hc}
\rm
For $\lambda \in \Lambda$ let
${\tilde T}_\lambda$ denote the
image of 
${\check T}_\lambda$ under 
the homomorphism $\varphi$ from
Lemma
\ref{lem:cancheck}.
\end{definition}

\begin{proposition}
\label{lem:tildehc}
The $\F$-vector space $\tilde T$ decomposes as
\begin{eqnarray*}
{\tilde T} = \sum_{\lambda \in \Lambda} 
{\tilde T}_\lambda
\qquad \qquad  (\mbox{\rm direct sum}).
\end{eqnarray*}
\end{proposition}
\noindent {\it Proof:}
Recall the map
$\varphi:
{\check T} \to 
{\tilde T} $
from
Lemma
\ref{lem:cancheck}.
To get 
${\tilde T} = \sum_{\lambda \in \Lambda} 
{\tilde T}_\lambda$,
apply $\varphi$ to each
side of
(\ref{eq:tcheckds}) and evaluate the
result using
Definition
\ref{def:hc}
and the surjectivity
of
$\varphi $.
To see that the sum
$\sum_{\lambda \in \Lambda} 
{\tilde T}_\lambda $
is direct, we pick any
$\lambda \in \Lambda$ and show 
that
${\tilde T}_\lambda $ has zero intersection with
$\sum_{\chi \in \Lambda\backslash \lambda} 
{\tilde T}_\chi$.
To this end 
we fix $u$ in the intersection and show
$u=0$.
By Definition 
\ref{def:hc} 
and since $u \in {\tilde T}_\lambda$,
there exists
$v \in 
 {\check T}_\lambda $
such that
$\varphi(v)=u$.
By Definition 
\ref{def:hc}
and since 
$ u \in \sum_{\chi \in \Lambda\backslash \lambda} 
{\tilde T}_\chi $,
there exists
$v' \in 
\sum_{\chi \in \Lambda\backslash \lambda} 
{\check T}_\chi $ such that
$\varphi(v')=u$.
Observe that $\varphi(v-v')=0$
so $v -v'\in R$.
By Lemma
\ref{def:relvector2}
there exists
$r \in R_\lambda$ and
$r' \in 
\sum_{\chi \in \Lambda\backslash \lambda} 
R_\chi $ such that
$v-v'=r-r'$.
Observe that
$v-r=v'-r'$.
We have
$v-r \in 
{\check T}_\lambda$
by construction and
the last sentence of
Definition
\ref{def:relatorsub}.
Similarly
$v'-r' \in
\sum_{\chi \in \Lambda\backslash \lambda} {\check T}_\chi$.
Now $v=r$ and $v'=r'$ in view of
Proposition
\ref{thm:dec}.
In the equation
$v=r$ we apply $\varphi$ to each side and get
$u=0$, as desired.
We have shown that the sum
$\sum_{\lambda \in \Lambda} 
{\tilde T}_\lambda $
is direct.
\hfill $\Box$ \\

\begin{definition}
\rm
For $\lambda \in \Lambda $ we call
${\tilde T}_\lambda$ the
{\it $\lambda$-homogeneous component} of 
$\tilde T$.
Elements of 
${\tilde T}_\lambda$ are said to be
{\it $\lambda$-homogeneous}.
An element of
${\tilde T}$ is called
{\it homogeneous} whenever
it is 
{\it $\lambda$-homogeneous}
for some $\lambda \in \Lambda$.
\end{definition}


\begin{proposition}
{\rm \cite[Theorem~8.1]{NT:onshape}}
\label{thm:zz} 
For all $\lambda \in \Lambda$ 
the map $\varphi$ sends $Z_\lambda$ onto
${\tilde T}_\lambda$.
\end{proposition}
\noindent {\it Proof:}
In \cite[Theorem~8.1]{NT:onshape} it is proved
that for an integer $n\geq 1$ and idempotent generators
$y,z$ of $T$ the following sets have the
same span:
\begin{enumerate}
\item[\rm (i)] The words of length $n$ in $T$ that begin
with $y$ and end with $z$.
\item[\rm (ii)] The zigzag words of length $n$ in $T$
that begin with $y$ and end with $z$.
\end{enumerate}
In that proof the relations
  (\ref{eq:sumeicop}) were never used; consequently
the verbatim proof applies to $\tilde T$ as well,
provided that we interpret things using
Note
\ref{ameaning}. The result follows.
\hfill $\Box$ \\

\begin{lemma}
\label{cor:zplusr}
For $\lambda \in \Lambda$,
\begin{eqnarray}
\label{eq:dirsum}
{\check T}_\lambda
=
R_\lambda + Z_\lambda.
\end{eqnarray}
\end{lemma}
\noindent {\it Proof:} 
The space
${\check T}_\lambda$
contains
$R_\lambda$
 by
Definition
\ref{def:relatorsub},
and it contains
$Z_\lambda $ by
Definition
\ref{def:zlambda}.
Consider the
 map
$\varphi : {\check T} \to {\tilde T}$
from Lemma
\ref{lem:cancheck}.
By 
Definition
\ref{def:hc}
 ${\tilde T}_\lambda$
is the image of
${\check T}_\lambda$ under 
$\varphi $. 
By Corollary
\ref{cor:rint}
$R_\lambda $ is the
kernel of 
$\varphi $ on 
${\check T}_\lambda$.
By Proposition
\ref{thm:zz}  $\varphi$ sends
$Z_\lambda$ onto 
${\tilde T}_\lambda$.
The result follows.
\hfill $\Box$ \\

\noindent We conjecture that the
sum
(\ref{eq:dirsum}) is direct for all
$\lambda \in \Lambda$.
For our present purpose the following
weaker result will suffice.
As part of our proof of
Theorem \ref{conj:mainp}  
we will show that
the sum (\ref{eq:dirsum}) is direct for all
$\lambda \in \Lambda_0$.
We will say more about this in the next section.
For the rest of this section we discuss
some aspects of (\ref{eq:dirsum}) that apply to
all $\lambda \in \Lambda$.

\begin{lemma}
\label{lem:dimbound}
The following hold for 
all $\lambda \in \Lambda $.
\begin{enumerate}
\item[\rm (i)]
${\rm dim} \, {\check T}_\lambda
=
{\rm dim} \, R_\lambda 
+
{\rm dim} \,Z_\lambda
-
{\rm dim} \,(R_\lambda\cap Z_\lambda)$.
\item[\rm (ii)]
${\rm dim} \,{\tilde T}_\lambda =
{\rm dim} \, Z_\lambda
-
{\rm dim} \,(R_\lambda\cap Z_\lambda)$.
\end{enumerate}
\end{lemma}
\noindent {\it Proof:}
(i) By Lemma
\ref{cor:zplusr} and elementary linear algebra.
\\
\noindent (ii) The action of $\varphi$  on 
$Z_\lambda$ is onto
${\tilde T}_\lambda$ and has kernel
$R_\lambda\cap Z_\lambda $.
\hfill $\Box$ \\

\begin{corollary}
\label{lem:dimbound2}
The following hold for 
all $\lambda \in \Lambda $.
\begin{enumerate}
\item[\rm (i)]
${\rm dim} \, {\check T}_\lambda
\leq 
{\rm dim} \, R_\lambda 
+
{\rm dim} \,Z_\lambda$.
\item[\rm (ii)]
${\rm dim} \, {\tilde T}_\lambda 
\leq 
{\rm dim} \,Z_\lambda$.
\end{enumerate}
\end{corollary}
\noindent {\it Proof:}
Immediate from Lemma
\ref{lem:dimbound}.
\hfill $\Box$ \\

\noindent
Pick $\lambda \in \Lambda $ and consider
when is the sum
(\ref{eq:dirsum}) direct.
Recall the sequence 
$p=(\lbrace \theta_i\rbrace_{i=0}^d;
\lbrace \theta^*_i\rbrace_{i=0}^d)$ 
from the first paragraph of this section.
Since $R_\lambda$
depends on 
$p$, it is conceivable
that the 
sum (\ref{eq:dirsum}) is direct for some
values of $p$ but not others. It is also 
conceivable that the field $\F$ matters.
The following definition will facilitate our
discussion of these issues.

\begin{definition}
\label{def:fav}
\rm 
For $\lambda \in \Lambda$, we say that
$\lambda $ is 
{\it $(p,\F)$-direct}
whenever
the sum 
(\ref{eq:dirsum}) is direct.
\end{definition}

\begin{lemma}
\label{lem:desirablemeaning}
For $\lambda \in \Lambda$ the following
are equivalent:
\begin{enumerate}
\item[\rm (i)] $R_\lambda \cap Z_\lambda=0$;
\item[\rm (ii)] $\lambda$ is $(p,\F)$-direct;
\item[\rm (iii)] 
equality holds in Corollary \ref{lem:dimbound2}(i);
\item[\rm (iv)] 
equality holds in Corollary \ref{lem:dimbound2}(ii);
\item[\rm (v)] the restriction of $\varphi $ to
$ Z_\lambda $ is injective.
\end{enumerate}
\end{lemma}
\noindent {\it Proof:}
${\rm (i)} \Leftrightarrow {\rm (ii)}$
By Definition
\ref{def:fav}.
\\
\noindent
${\rm (i)} \Leftrightarrow {\rm (iii)}$
By Lemma
\ref{lem:dimbound}(i).
\\
\noindent
${\rm (i)} \Leftrightarrow {\rm (iv)}$
By Lemma
\ref{lem:dimbound}(ii).
\\
\noindent
${\rm (i)} \Leftrightarrow {\rm (v)}$
The restriction of $\varphi$ to $Z_\lambda$ has
kernel
 $R_\lambda \cap Z_\lambda$.
\hfill $\Box$ \\

\begin{proposition}
\label{prop:sometimes}
Assume $\F$ is infinite
and pick an integer $d\geq 3$.
Suppose we are given
a type $\lambda \in \Lambda$ that is
$(p,\F)$-direct
for all
sequences $p \in {\rm Rac}(d,\F)$.
Then
$\lambda $ is
$(p,\F)$-direct
for all sequences
$p \in {\rm Feas}(d,\F)$.
\end{proposition}
\noindent {\it Proof:}
For notational convenience abbreviate
$d_\lambda = 
{\rm dim} \, {\check T}_\lambda
-
{\rm dim} \, Z_\lambda
$.
Recall the relators 
 of $\check T$ from
Definition
\ref{def:relator}.
In the definition of
a relator
an element $p \in
{\rm Feas}(d,\F)$ is involved,
so that relator 
can be viewed as a function
of $p$.
We adopt this point of view throughout the proof. 
 Let ${\mathcal R}_\lambda$
denote the set of all $\lambda$-homogeneous
relators in $\check T$.
By Definition
\ref{def:relatorsub}
we have
$R_\lambda(p)={\rm Span}\lbrace 
\varrho(p) | 
\varrho \in {\mathcal R}_\lambda \rbrace$
for 
all $p \in
{\rm Feas}(d,\F)$.
We assume that there exists
$p' \in {\rm Feas}(d,\F)$
such that
 $\lambda$ is not $(p',\F)$-direct, and get a contradiction.
By Lemma
\ref{lem:desirablemeaning}(ii),(iii)
we have
${\rm dim} \, R_\lambda(p') \geq d_\lambda +1$.
By our above comments
$R_\lambda(p')={\rm Span}\lbrace 
\varrho(p') | 
\varrho \in {\mathcal R}_\lambda \rbrace$.
Therefore there exists a 
 subset 
$H \subseteq {\mathcal R}_\lambda$
such that (i) $H$ has cardinality 
$d_\lambda+1$; and (ii)
the set $\lbrace \varrho(p')\rbrace_{\varrho \in H}$
is linearly independent.
Pick any
$p \in {\rm Feas}(d,\F)$.
Recall by Definition
\ref{def:tlambda}
that ${\check T}_\lambda$ has a basis
consisting of the words of type $\lambda$.
For
$\varrho \in H$
write
$\varrho(p)$
as a linear
combination of these words,
and let $M=M(p)$ denote the corresponding
coefficient matrix.
The rows of $M$ are indexed by the
words of type $\lambda$,
and the columns of
$M$ are indexed by $H$.
Each entry of $M$ is a power of
some
$\theta_i$ or $\theta^*_i$,
where
   $p=(\lbrace \theta_i\rbrace_{i=0}^d;
\lbrace \theta^*_i\rbrace_{i=0}^d)$.
The matrix $M(p')$ has full rank
$d_\lambda+1$
since $\lbrace \varrho(p')\rbrace_{\varrho \in H}$
are linearly independent.
Therefore there exists a
set $L$ consisting of words of type
$\lambda$ such that 
(i) $L$ has cardinality
$d_\lambda+1$; 
and (ii) 
the rows of $M(p')$ indexed by $L$ are
linearly independent.
For
$p \in {\rm Feas}(d,\F)$ let
$N=N(p)$ denote 
the submatrix of $M(p)$ obtained
by deleting all rows
not indexed by
$L$.
By construction $N$ is 
$(d_\lambda+1)\times (d_\lambda+1)$,
and its determinant
 is a polynomial
in $p$ that has all coefficients
in $\F$.
Denote this polynomial by $h$.
By construction $N(p')$ is nonsingular
so $h(p')\not=0$.
We will obtain a contradiction by showing 
that $h(p')=0$. To this end
we will show that $h(p)=0$ for all
$p \in 
{\rm Rac}(d,\F)$, 
and invoke Proposition
\ref{lem:bootstrap}.
Pick any 
$p \in 
{\rm Rac}(d,\F)$.
By assumption $\lambda $ is $(p,\F)$-direct.
So 
 ${\rm dim} \, R_\lambda(p) = d_\lambda$
in view of Lemma
\ref{lem:desirablemeaning}(ii),(iii).
Now 
$N(p)$ is singular and hence
$h(p)=0$.
We have shown that $h(p)=0$ for
all $p \in 
{\rm Rac}(d,\F)$.
Now by Proposition
\ref{lem:bootstrap},  $h(p)=0$ for all 
$p \in {\rm Feas}(d,\F)$. In particular
$h(p')=0$, for a contradiction. 
The result follows.
\hfill $\Box$ \\

\section{The algebra
$\varepsilon^*_0 {\tilde T}\varepsilon^*_0$} 

Throughout this section
 fix an integer $d\geq 0$ and
  a sequence  $p=(\lbrace \theta_i\rbrace_{i=0}^d;
\lbrace \theta^*_i\rbrace_{i=0}^d)$ 
in
${\rm Feas}(d,\F)$.
Recall the algebra
${\check T}={\check T}(d,\F)$
from Definition
\ref{def:free} and
 the algebra  
${\tilde T}={\tilde T}(p,\F)$ from
 Definition
\ref{def:tt}.
Observe that $\varepsilon^*_0 {\tilde T}\varepsilon^*_0$ 
is an $\F$-algebra
with multiplicative identity 
 $\varepsilon^*_0$.

\medskip
\noindent Recall the map
$\varphi : {\check T} \to {\tilde T}$
from Lemma
\ref{lem:cancheck}.

\begin{lemma}
\label{lem:etazero}
The restriction of
$\varphi$ to 
$\epsilon^*_0 {\check T}\epsilon^*_0$
gives a surjective $\F$-algebra homomorphism
$\epsilon^*_0 {\check T}\epsilon^*_0 \to
\varepsilon^*_0 {\tilde T}\varepsilon^*_0$.
\end{lemma}
\noindent {\it Proof:}
By Lemma
\ref{lem:cancheck}
the map
$\varphi : {\check T} \to {\tilde T}$
is a surjective 
$\F$-algebra homomorphism
that sends
$\epsilon^*_0 \mapsto
\varepsilon^*_0$. The result follows.
\hfill $\Box$ \\

\begin{proposition}
\label{cor:tzerodec}
The $\F$-vector space 
$\varepsilon^*_0{\tilde T}\varepsilon^*_0$ decomposes as
\begin{eqnarray}
\label{eq:tildetzerodec}
\varepsilon^*_0{\tilde T}\varepsilon^*_0
=
\sum_{n=0}^\infty
{\tilde T}_{\lbrack n \rbrack}
\qquad \qquad  (\mbox{\rm direct sum}).
\end{eqnarray}
Moreover 
${\tilde T}_{\lbrack m \rbrack}
\cdot
{\tilde T}_{\lbrack n \rbrack}
\subseteq 
{\tilde T}_{\lbrack m+n \rbrack}
$
for all integers $m,n\geq 0$.
\end{proposition}
\noindent {\it Proof:}
To get
$\varepsilon^*_0{\tilde T}\varepsilon^*_0
=
\sum_{n=0}^\infty
{\tilde T}_{\lbrack n \rbrack}$,
 apply $\varphi$ to each side of
(\ref{eq:tzerodec}) and evaluate the
result using
Definition
\ref{def:hc}
and Lemma
\ref{lem:etazero}.
The sum 
$\sum_{n=0}^\infty
{\tilde T}_{\lbrack n \rbrack}$ is direct
by 
Proposition
\ref{lem:tildehc}.
The last assertion follows from
the last assertion of Proposition
\ref{lem:ete1}.
\hfill $\Box$ \\

\begin{lemma} 
\label{lem:tmucom}
The elements 
$\lbrace \varepsilon^*_0 \varepsilon_i \varepsilon^*_0\rbrace_{i=0}^d $
mutually commute.
\end{lemma}
\noindent {\it Proof:}
In 
\cite[Theorem~2.4]{nomtowards} it was proved
that 
$\lbrace e^*_0 e_i e^*_0\rbrace_{i=0}^d $
commute in $T$.
In that proof the relations
(\ref{eq:sumeicop}) were never used.
Consequently the verbatim proof applies
to the elements
$\lbrace \varepsilon^*_0 \varepsilon_i \varepsilon^*_0\rbrace_{i=0}^d $
of $\tilde T$,
provided that we interpret things using
Note
\ref{ameaning}.
\hfill $\Box$ \\

\noindent Let
$\lbrace x_i \rbrace_{i=0}^d $ denote mutually
commuting indeterminates. Let
$\F\lbrack x_0,\ldots, x_d\rbrack$ denote
the $\F$-algebra consisting of the
polynomials in 
$\lbrace x_i\rbrace_{i=0}^d$
that have all coefficients in $\F$.
We abbreviate $P=
\F\lbrack x_0,\ldots, x_d\rbrack$.

\begin{corollary}
\label{cor:surtilde}
There exists an $\F$-algebra 
homomorphism ${\tilde \nu} : P
\to 
\varepsilon^*_0{\tilde T}\varepsilon^*_0$ 
that sends $x_i \mapsto 
\varepsilon^*_0 \varepsilon_i \varepsilon^*_0 $ for $0 \leq i \leq d$.
\end{corollary}
\noindent {\it Proof:}  Immediate from
Lemma
\ref{lem:tmucom}.
\hfill $\Box$ \\

\noindent In Section 18  we
will show
that
the map $\tilde \nu$ from Corollary
\ref{cor:surtilde} is an isomorphism.
For the time being, our goal is
to show that 
$\tilde \nu$ is an isomorphism if and only if
$\lbrack n \rbrack $ is $(p,\F)$-direct 
 for all integers $n\geq 0$.

\medskip
\noindent As we discuss the algebra $P$
the following notation will be helpful.
We call $\lbrace x_i\rbrace_{i=0}^d$ the
{\it generators} for $P$.
For a generator $x_i$ of $P$ we call $i$ the {\it index} of
$x_i$. For a generator $y$ of $P$ let $\overline y$ denote
the index of $y$.
For an integer $n\geq 0$, by a {\it monomial of degree $n$}
in $P$ we mean an element $y_1y_2\cdots y_n$
such that 
$y_i$ is a generator of $P$ for $1 \leq i \leq n$.
For notational convenience 
we always order the factors such that
${\overline y}_{i-1} \geq {\overline y}_i$
for $2 \leq i \leq n$.
We interpret the monomial of degree 0 to be the identity
of $P$.
Observe that the $\F$-vector space
$P$ has a basis consisting of its monomials.
For $n\geq 0$ let $P_n$ denote the subspace of
 $P$ with a basis consisting of the
monomials of degree $n$. We have
\begin{eqnarray}
P = \sum_{n=0}^\infty P_n \qquad \qquad 
{\mbox{\rm (direct sum).}}
\label{eq:pdecomp}
\end{eqnarray}
Moreover $P_mP_n=P_{m+n}$ for all $m,n\geq 0$.
We call $P_n$ the {\it nth homogeneous component}
of $P$.

\begin{definition}
\label{def:ximap}
\rm
We define an $\F$-linear map
$\natural : P \to \epsilon^*_0{\check T}\epsilon^*_0$.
To do this we give the action
of $\natural$ on the monomial basis for $P$.
By definition $\natural$ sends
$1\mapsto \epsilon^*_0$.
For each nontrivial monomial
$y_1y_2\cdots y_n $ in $P$ the
image under $\natural$ is 
$\epsilon^*_0
y'_1
\epsilon^*_0
y'_2
\epsilon^*_0
\cdots
\epsilon^*_0
y'_n
\epsilon^*_0
$, where
$x'_i = \epsilon_i$ for $0 \leq i \leq d$.
\end{definition}

\noindent We caution the reader that
$\natural$
is not an algebra homomorphism in general.

\begin{lemma}
\label{lem:xipre}
For an integer $n\geq 0$ 
the map $\natural$ from Definition
\ref{def:ximap}
induces a bijection between
the following two sets:
\begin{enumerate}
\item[\rm (i)] the monomials in $P$ that have degree $n$;
\item[\rm (ii)] 
the zigzag words in $\check T$ of type $\lbrack n \rbrack$.
\end{enumerate}
\end{lemma}
\noindent {\it Proof:}
Compare Proposition
\ref{prop:zzn}
and Definition
\ref{def:ximap}.
\hfill $\Box$ \\

\begin{lemma}
\label{lem:bijbij}
The map
$\natural$ from
Definition
\ref{def:ximap}
is an injection. 
For $n\geq 0$ the image
of $P_n$ under 
$\natural$ is equal to
$Z_{\lbrack n \rbrack}$.
\end{lemma}
\noindent {\it Proof:}
The monomials in $P$ form a basis for $P$.
The zigzag words in $\check T$ are linearly independent.
By these comments and Lemma
\ref{lem:xipre}
the map $\natural$ is injective.
To get the last assertion, note that in
Lemma
\ref{lem:xipre}
the set (i) is a basis for $P_n$ and
the set (ii) is a basis for $Z_{\lbrack n \rbrack}$.
\hfill $\Box$ \\

\begin{lemma} 
\label{lem:diagcomm}
Let $\varphi'$ denote the restriction of
$\varphi$ to
$\epsilon^*_0 {\check T}\epsilon^*_0$.
Then the following diagram commutes:

\[
\begin{CD}
P  @>{\rm id}>>
         P \\
	  @V\natural VV                     @VV\tilde \nu V \\
	\epsilon^*_0{\check T}\epsilon^*_0 
	  @>>\varphi'>
	           \varepsilon^*_0{\tilde T}\varepsilon^*_0 
		    \end{CD}
		    \]

\end{lemma}
\noindent {\it Proof:}
The map $\natural$ is  from
Definition
\ref{def:ximap} and the map
$\tilde \nu$ is from
Corollary
\ref{cor:surtilde}. 
The map
$\varphi$ is from
Lemma
\ref{lem:cancheck}.
The result
is a
routine consequence of how these maps
are defined.
\hfill $\Box$ \\

\begin{lemma}
\label{lem:hominj}
For an integer $n\geq 0$ 
the image of $P_n$ under $\tilde \nu$
is ${\tilde T}_{\lbrack n \rbrack}$.
\end{lemma}
\noindent {\it Proof:}
By
Lemma \ref{lem:diagcomm}
the composition $\varphi' \circ \natural$ equals
 ${\tilde \nu}$.
By Lemma
\ref{lem:bijbij} the image of 
$P_n$ under 
$\natural$ is 
$Z_{\lbrack n \rbrack}$.
By
Lemma
\ref{thm:zz}
the space ${\tilde T}_{\lbrack n \rbrack}$
is the image of
$Z_{\lbrack n \rbrack}$ under 
$\varphi $ and hence 
$\varphi' $.
The result follows.
\hfill $\Box$ \\

\begin{lemma} 
\label{lem:tmu}
The map
 ${\tilde \nu}$
from
Corollary
\ref{cor:surtilde}
 is surjective.
\end{lemma}
\noindent {\it Proof:}
In the equation
(\ref{eq:pdecomp}) apply
$\tilde \nu$ to each side, and evaluate the result
using
Lemma
\ref{lem:hominj}
and then
Proposition
\ref{cor:tzerodec}.
\hfill $\Box$ \\

\begin{lemma}
\label{prop:equivpre}
For an integer $n\geq 0$ the following are equivalent:
\begin{enumerate}
\item[\rm (i)]
the restriction ${\tilde \nu}$ to $P_n$ is injective;
\item[\rm (ii)]
the type $\lbrack n \rbrack $ is $(p,\F)$-direct
in the sense of Definition
\ref{def:fav}.
\end{enumerate}
\end{lemma}
\noindent {\it Proof:}
Consider the commuting diagram in Lemma
\ref{lem:diagcomm}.
By Lemma
\ref{lem:bijbij}
the map $\natural$ is an injection that
sends $P_n$ onto $Z_{\lbrack n \rbrack}$.
Therefore the restriction of $\tilde \nu$ to $P_n$ is
injective if and only if the restriction of
$\varphi$ to $Z_{\lbrack n\rbrack}$ is injective.
By Lemma
\ref{lem:desirablemeaning}
the restriction of $\varphi$ to 
$Z_{\lbrack n \rbrack}$ is injective if and only
if $\lbrack n \rbrack$ is $(p,\F)$-direct.
The result follows.
\hfill $\Box$ \\

\begin{lemma}
\label{lem:kerintersect}
The kernel of
${\tilde \nu}$
decomposes as follows:
\begin{eqnarray}
\label{eq:sumsum}
{\rm ker}({\tilde \nu})
=
\sum_{n=0}^\infty
\bigl({\rm ker}({\tilde \nu}) \cap
P_n \bigr).
\end{eqnarray}
\end{lemma}
\noindent {\it Proof:} 
The inclusion
$\supseteq$ is clear, so consider the
inclusion
$\subseteq$.
Pick $h \in 
{\rm ker}({\tilde \nu})$.
By 
(\ref{eq:pdecomp})
there exists an integer $m\geq 0$ 
and a sequence
$\lbrace h_n\rbrace_{n=0}^m$
 such that $h_n \in
P_n$
for $0 \leq n\leq m$ and
$h=\sum_{n=0}^m h_n$.
In this equation we apply 
 $\tilde \nu$ to each term
and get
$0 = \sum_{n=0}^m {\tilde \nu}(h_n)$.
By Lemma \ref{lem:hominj} we have
${\tilde \nu}(h_n) \in {\tilde T}_{\lbrack n \rbrack}$
for $0 \leq n \leq m$.
By these comments and
(\ref{eq:tildetzerodec})
 we obtain
${\tilde \nu}(h_n)=0$ for  $ 0 \leq n \leq m$.
So for $0 \leq n \leq m$ the
polynomial
$h_n$ is contained in the $n$-summand on
the right in
(\ref{eq:sumsum}).
Therefore $h$  is contained in the sum on
the right in (\ref{eq:sumsum}).
We have verified the inclusion
$\subseteq$ and the result follows.
\hfill $\Box$ \\

\begin{proposition}
\label{prop:equiv}
The following are equivalent:
\begin{enumerate}
\item[\rm (i)]
the map $\tilde \nu$ from Corollary
\ref{cor:surtilde} is an isomorphism;
\item[\rm (ii)]
for all integers $n\geq 0$
the type $\lbrack n \rbrack $ is $(p,\F)$-direct
in the sense of
Definition
\ref{def:fav}.
\end{enumerate}
\end{proposition}
\noindent {\it Proof:}
The map $\tilde \nu$ is surjective by Lemma
\ref{lem:tmu}, so
$\tilde \nu$ is an isomorphism
if and only if
$\tilde \nu$ is injective.
By Lemma
\ref{lem:kerintersect}
the map $\tilde \nu$ is injective
if and only if its restriction to
$P_n$ is
injective for all $n\geq 0$.
The result follows from these comments and
Lemma
\ref{prop:equivpre}.
\hfill $\Box$ \\

\section{A central element of ${\tilde T}$}

\noindent Recall that the algebra $T$
from Definition
\ref{def:twitha} is defined using
relations
(\ref{eq:eiejcop})--(\ref{eq:eiaskej}).
So far we have investigated all these relations
except 
(\ref{eq:sumeicop}). We now prepare to bring in
the relations (\ref{eq:sumeicop}).

\medskip
\noindent
Throughout this section we fix an integer $d\geq 0$ and
 a  sequence  $p=(\lbrace \theta_i\rbrace_{i=0}^d;
\lbrace \theta^*_i\rbrace_{i=0}^d)$ 
in ${\rm Feas}(d,\F)$. 
Recall the algebra ${\tilde T}={\tilde T}(p,\F)$ 
from Definition
\ref{def:tt}.

\begin{definition}
\label{def:del}
\rm Define $\Delta \in {\tilde T}$ and
 $\Delta^* \in {\tilde T}$ by
\begin{eqnarray*}
\Delta = 1-\sum_{i=0}^d \varepsilon_i,
\qquad \qquad 
\Delta^*= 1-\sum_{i=0}^d \varepsilon^*_i.
\end{eqnarray*}
\end{definition}

\noindent 
The elements 
 $\Delta$ and $\Delta^*$
are nonzero by Proposition
\ref{lem:tildehc}.

\begin{lemma}
\label{lem:del1}
We have 
$\Delta^2 =\Delta $ and
$\Delta^{*2} =\Delta^*$. Moreover
\begin{eqnarray}
\label{eq:eprod1}
\varepsilon_i \Delta= \Delta \varepsilon_i = 0,
\qquad \qquad
\varepsilon^*_i \Delta^*= \Delta^* \varepsilon^*_i = 0
\qquad \qquad
0 \leq i \leq d.
\end{eqnarray}
\end{lemma}
\noindent {\it Proof:}
Line
(\ref{eq:eprod1}) follows
from
(\ref{eq:r1}) and
  Definition
\ref{def:del}.
To obtain 
$\Delta^2 =\Delta $,
observe that 
$\Delta (1-\Delta)=\sum_{i=0}^d \Delta \varepsilon_i=0$.
The equation
$\Delta^{*2} =\Delta^* $ 
is similarly obtained.
\hfill $\Box$ \\

\begin{lemma}
\label{lem:del1a}
For $0 \leq i,j\leq d$ with $i\not=j$,
\begin{eqnarray*}
\varepsilon^*_i\Delta \varepsilon^*_j = 0,
\qquad \qquad 
\varepsilon_i\Delta^* \varepsilon_j = 0.
\end{eqnarray*}
\end{lemma}
\noindent {\it Proof:}
By Definition 
\ref{def:del}
and since $\vare^*_i \vare^*_j=0$, we find
$\vare^*_i\Delta \vare^*_j=-\sum_{\ell=0}^d \vare^*_i \vare_\ell \vare^*_j$.
Setting $k=0$ in the equation on the left in 
(\ref{eq:r2}) we find
$0=\sum_{\ell=0}^d \vare^*_i \vare_\ell \vare^*_j$.
Therefore
$\vare^*_i\Delta \vare^*_j = 0$. 
The equation
$\vare_i\Delta^* \vare_j = 0$ is similarly obtained.
\hfill $\Box$ \\

\begin{definition}
\label{def:psi}
\rm
Define $\psi \in {\tilde T}$ by
\begin{eqnarray}
\psi = (\Delta-\Delta^*)^2,
\label{eq:psi}
\end{eqnarray}
where $\Delta$, $\Delta^*$ are from Definition
\ref{def:del}.
\end{definition}

\noindent An element of an algebra is called
{\it central} whenever it commutes with
everything in the algebra.
Our next goal is to show that $\psi$ is central.

\begin{lemma}
\label{lem:centpre}
The element $\psi$ coincides with
each of the following:
\begin{enumerate}
\item[\rm (i)] $\Delta+\Delta^*-\Delta \Delta^* -\Delta^*\Delta$;
\item[\rm (ii)] $\Delta (\vare^*_0 + \vare^*_1+ \cdots + \vare^*_d)+
 \Delta^* (\vare_0 + \vare_1+ \cdots + \vare_d)$;
\item[\rm (iii)] $(\vare^*_0 + \vare^*_1+ \cdots + \vare^*_d)\Delta+
 (\vare_0 + \vare_1+ \cdots + \vare_d)\Delta^*$.
\end{enumerate}
\end{lemma}
\noindent {\it Proof:}
(i)
Multiply out the right-hand side of 
(\ref{eq:psi}), and simplify the result using
$\Delta^2=\Delta$ and
$\Delta^{*2}=\Delta^*$.
\\
\noindent (ii) In the given expression eliminate 
$\vare_0+\vare_1 + \cdots + \vare_d$ and
$\vare^*_0+\vare^*_1 + \cdots + \vare^*_d$
using Definition
\ref{def:del}, and compare the result with
(i) above.
\\
\noindent (iii) Similar to the proof of (ii) above.
\hfill $\Box$ \\

\begin{lemma}
\label{lem:del1half}
The following 
hold for $0 \leq i \leq d$.
\begin{enumerate}
\item[\rm (i)] Each of $\vare_i\psi$, $\psi \vare_i$ is equal to
$\vare_i \Delta^* \vare_i$.
\item[\rm (ii)] Each of $\vare^*_i \psi$, $\psi \vare^*_i$ is equal to
$\vare^*_i \Delta \vare^*_i$.
\end{enumerate}
\end{lemma}
\noindent {\it Proof:}
(i) Evaluating 
$\vare_i\psi$
using
Lemma
\ref{lem:centpre}(ii) and $\vare_i\Delta = 0$,
we obtain
$\vare_i \psi =
\vare_i \Delta^* (\vare_0 + \vare_1 + \cdots + \vare_d)$.
By
Lemma \ref{lem:del1a},
for $0 \leq j \leq d$
we have
$\vare_i \Delta^* \vare_j = 0 $ provided $i \not= j$.
By these comments
$\vare_i \psi = \vare_i \Delta^* \vare_i$. Using
Lemma
\ref{lem:centpre}(iii) we similarly find
$\psi \vare_i = \vare_i \Delta^* \vare_i$. 
\\
\noindent (ii) Similar to the proof of (i) above.
\hfill $\Box$ \\

\begin{corollary}
\label{lem:del2}
The element $\psi$ 
is central in ${\tilde T}$.
\end{corollary}
\noindent {\it Proof:}
The elements $\lbrace \vare_i\rbrace_{i=0}^d$,
$\lbrace \vare^*_i\rbrace_{i=0}^d$ together generate
${\tilde T}$, and each of these 
elements commutes with $\psi$ by Lemma
\ref{lem:del1half}.
\hfill $\Box$ \\

\noindent For later use we summarize
Lemma \ref{lem:del1a} and
Lemma
\ref{lem:del1half}.

\begin{lemma}
\label{lem:psie}
For $0 \leq i,j\leq d$ we have
\begin{eqnarray*}
&&
\vare_i \Delta^* \vare_j = \delta_{i,j} \psi \vare_i,
\qquad \qquad
\vare^*_i \Delta \vare^*_j = \delta_{i,j} \psi \vare^*_i.
\end{eqnarray*}
\end{lemma}

\section{The 
 homomorphism 
$\pi: {\tilde T}\to T$}

\noindent
Throughout this section 
 fix an integer $d\geq 0$ and
 a sequence  $p=(\lbrace \theta_i\rbrace_{i=0}^d;
\lbrace \theta^*_i\rbrace_{i=0}^d)$ 
in
${\rm Feas}(d,\F)$.
Recall the algebras
$ T= T(p,\F)$ 
from 
Definition
\ref{def:twitha}
and 
${\tilde T}={\tilde T}(p,\F)$ 
from Definition
\ref{def:tt}.

\begin{definition}
\label{def:j}
\rm 
Let $J$ denote the two-sided ideal of ${\tilde T}$ 
generated by the elements $\Delta, \Delta^*$ from
Definition
\ref{def:del}.
\end{definition}

\begin{lemma}
\label{lem:ideal}
There exists a surjective $\F$-algebra homomorphism
$\pi: {\tilde T} \to T$ that sends
$\varepsilon_i \mapsto e_i$ and
$\varepsilon^*_i \mapsto e^*_i$ for $0 \leq i \leq d$.
The kernel of $\pi$ coincides with the ideal $J$.
\end{lemma}
\noindent {\it Proof:}
Compare the defining relations
for 
${\tilde T}$ and $T$.
\hfill $\Box$ \\

\noindent
We will be discussing the action 
of 
$\pi$ on  $\varepsilon^*_0{\tilde T}\varepsilon^*_0$.

\begin{lemma}
\label{lem:pisurj}
\label{def:piprime}
The restriction of
$\pi$
to 
$\vare^*_0{\tilde T}\vare^*_0$
gives a surjective $\F$-algebra homomorphism
$\vare^*_0{\tilde T}\vare^*_0 \to
e^*_0 Te^*_0$.
\end{lemma}
\noindent {\it Proof:}
The map
$\pi$
is a surjective $\F$-algebra
homomorphism that sends
$\vare^*_0 \mapsto 
 e^*_0$.
\hfill $\Box$ \\

\begin{proposition}
\label{thm:idealinter}
The following are equal:
\begin{enumerate}
\item[\rm (i)] 
the kernel of 
$\pi$ on  $\varepsilon^*_0{\tilde T}\varepsilon^*_0$;
\item[\rm (ii)] 
the intersection of $J$ and $\vare^*_0 {\tilde T}\vare^*_0$;
\item[\rm (iii)] 
 $\vare^*_0 J\vare^*_0$;
\item[\rm (iv)]
the ideal of $\vare^*_0{\tilde T}\vare^*_0$ generated
by $\vare^*_0\Delta \vare^*_0$.
\end{enumerate}
\end{proposition}
\noindent {\it Proof:}
The spaces (i), (ii) are equal by
 the last assertion of Lemma
\ref{lem:ideal}.
We now show that the spaces
(ii)--(iv) are equal.
Let $J'$ denote the ideal of
 $\vare^*_0{\tilde T}\vare^*_0$ generated
by $\vare^*_0\Delta \vare^*_0$.
\\
\noindent 
$J\cap 
\vare^*_0 {\tilde T} \vare^*_0 \subseteq 
\vare^*_0 J \vare^*_0 $:
For $u \in 
\vare^*_0 {\tilde T} \vare^*_0 $
we have $u = 
\vare^*_0 u \vare^*_0$ since
$\vare^{*2}_0 =
\vare^*_0$.
\\
\noindent 
$\vare^*_0 J \vare^*_0 \subseteq J'$:
Let $J_1$ (resp. $J_2$)
denote the two-sided ideal of $\tilde T$
generated by $\Delta $ (resp. $\Delta^*$).
By construction $J=J_1+J_2$, so
$\vare^*_0J\vare^*_0=\vare^*_0J_1\vare^*_0+\vare^*_0J_2 \vare^*_0$.
We now show that 
$\vare^*_0 J_1 \vare^*_0
\subseteq J'$.
The space $ \vare^*_0 J_1 \vare^*_0$
is spanned by elements of the
form 
$u \Delta v$ where
$u$ (resp. $v$)
is a nontrivial word in $\tilde T$
that begins with $\vare^*_0$
(resp. ends with $\vare^*_0$).
We show that such an element
$u \Delta v$ is contained in
$J'$.
Suppose for the moment that
$u$ ends with a nonstarred idempotent generator, which we 
denote by $\vare_i$. Then $u\Delta=0$ since
$\vare_i\Delta=0$.
Therefore we may assume
that $u$ ends with a starred idempotent 
generator, which we denote $\vare^*_i$.
Note that $u=u\vare^*_i$ since $\vare^{*2}_i=\vare^*_i$.
Suppose for the moment that
$v$ begins with a
nonstarred idempotent generator, which we  denote by $\vare_j$.
 Then $\Delta v=0$ since
$\Delta \vare_j=0$.
Therefore we may assume
that $v$ begins with a starred idempotent generator, which we denote
by $\vare^*_j$.
Note that $v=\vare^*_jv$ since
 $\vare^{*2}_j=\vare^*_j$.
We may now argue
\begin{eqnarray*}
u\Delta v 
&=&
u\vare^*_i\Delta \vare^*_j v 
\\
&=&
\delta_{i,j}u \psi \vare^*_i v
\;
\;\quad \qquad \qquad {\mbox{(by
Lemma \ref{lem:psie})}}
\\
&=&
\delta_{i,j}\psi u  \vare^*_i v 
\;\quad \qquad \qquad {\mbox{ (by
Corollary \ref{lem:del2})}}
\\
&=&
\delta_{i,j}\psi u v.
\end{eqnarray*}
Since $u$ begins with $\vare^*_0$ and
$\vare^{*2}_0=\vare^*_0$ we find
$u=\vare^*_0u$. 
Also $\psi \vare^*_0=
 \vare^*_0\Delta \vare^*_0$ by
Lemma \ref{lem:psie}.
Therefore
$\psi uv=
\vare^*_0\Delta \vare^*_0 uv$.
Since $u$ begins with $\vare^*_0$ and 
$v$ ends with $\vare^*_0$ we find
$uv \in \vare^*_0{\tilde T} \vare^*_0$.
Consequently $J'$ contains
$\psi uv$ 
and hence
$u\Delta v$.
We have shown
$\vare^*_0 J_1\vare^*_0
\subseteq J'$.
One similarly shows that
$\vare^*_0 J_2 \vare^*_0
\subseteq J'$.
\\
\noindent
$J' \subseteq J\cap 
\vare^*_0 {\tilde T} \vare^*_0$:
Observe that $J' \subseteq 
J$ since
$\Delta \in J$,
and $J' \subseteq 
\vare^*_0{\tilde T} \vare^*_0$
by construction.
\hfill $\Box$ \\

\section{Comparing $\vare^*_0{\tilde T}\vare^*_0$ and
$e^*_0 Te^*_0$}

\noindent 
Throughout this section 
fix an integer $d\geq 0$ and
 a sequence  $p=(\lbrace \theta_i\rbrace_{i=0}^d;
\lbrace \theta^*_i\rbrace_{i=0}^d)$ 
in
${\rm Feas}(d,\F)$.
Recall the algebras
$ T= T(p,\F)$ 
from 
Definition
\ref{def:twitha}
and 
${\tilde T}={\tilde T}(p,\F)$ 
from Definition
\ref{def:tt}.
We will be comparing the map
$\mu:\F\lbrack x_1,\ldots, x_d\rbrack \to e^*_0Te^*_0$
from
Corollary \ref{cor:18a},
and the map
${\tilde \nu}:P \to 
\vare^*_0 {\tilde T}\vare^*_0$ from Corollary
\ref{cor:surtilde}. We will show that
 $\mu$ is an isomorphism if and only 
if $\tilde \nu$ is an isomorphism.

\medskip
\noindent  In order to
compare $\mu$ and $\tilde \nu$ it is helpful
to introduce the following map.

\begin{definition}
\label{def:nu}
\rm
Define an 
$\F$-algebra homomorphism
$\nu:
\F\lbrack x_1,  \ldots, x_d\rbrack \to 
e^*_0 Te^*_0$  that sends $x_i \to e^*_0e_ie^*_0$ for
$1 \leq i \leq d$.
\end{definition}

\noindent Our next goal is to compare
$\mu$ and $\nu$.
After that, we will compare
$\nu$ and $\tilde \nu$.

\begin{definition}
\label{def:phi}
\rm
Define an 
$\F$-algebra homomorphism
$\phi:
\F\lbrack x_1,  \ldots, x_d\rbrack \to 
\F\lbrack x_1,  \ldots, x_d\rbrack $
that sends $x_i\to
\sum_{j=1}^d \tau_i(\theta_j)x_j$
for
$1 \leq i \leq d$.
(The $\tau_i$ are from Definition
\ref{def:taueta}).
\end{definition}

\begin{lemma}
\label{lem:phiiso}
The map $\phi$ from Definition
\ref{def:phi} is an isomorphism.
\end{lemma}
\noindent {\it Proof:} 
Consider the $d \times d$ matrix
that has $(i,j)$-entry
$\tau_i(\theta_j)$ for $1 \leq i,j\leq d$.
This matrix is upper triangular and has all diagonal
entries nonzero.
Therefore the matrix is invertible.
The result follows.
\hfill $\Box$ \\

\begin{lemma}
\label{lem:diag1}
The following diagram commutes:

\[
\begin{CD}
\F\lbrack x_1,\ldots, x_d\rbrack  @>\phi>>
         \F\lbrack x_1,\ldots,x_d\rbrack \\
	  @V\mu VV                     @VV\nu V \\
	e^*_0 Te^*_0 
	  @>>{\rm id}>
	           e^*_0Te^*_0 
		    \end{CD}
		    \]

\end{lemma}

\noindent {\it Proof:} 
For $1 \leq i \leq d$ we chase $x_i$ around
the diagram.
The image of $x_i$ under
the composition
$\nu \circ \phi $ 
is
$\sum_{j=1}^d \tau_i(\theta_j)e^*_0e_je^*_0$.
The image of $x_i$ under 
$\mu$ is
$e^*_0\tau_i(a)e^*_0$,
and this is equal to 
$\sum_{j=0}^d \tau_i(\theta_j)e^*_0e_je^*_0$.
In this sum the $j=0$ summand is zero; indeed
$\tau_i(\theta_0)=0$ since $i\geq 1$.
Therefore
$x_i$ has the same image under
$\nu \circ \phi $ and
$\mu$.
The result follows.
\hfill $\Box$ \\

\begin{corollary}
\label{cor:nusurj}
The map
$\nu$ is surjective.
\end{corollary}
\noindent {\it Proof:} 
The map $\mu$ is surjective by Corollary
\ref{cor:18a}.
The result follows from this and
Lemma \ref{lem:diag1}.
\hfill $\Box$ \\

\begin{proposition}
\label{lem:isodual}
The following 
are equivalent:
\begin{enumerate}
\item[\rm (i)]
the map $\mu$
is an isomorphism; 
\item[\rm (ii)]
the map $\nu$
is an isomorphism.
\end{enumerate}
\end{proposition}
\noindent {\it Proof:} 
Combine
Lemma
\ref{lem:phiiso}
and 
Lemma
\ref{lem:diag1}.
\hfill $\Box$ \\

\noindent Our next goal is to compare
$\nu$ and $\tilde \nu$.

\begin{definition}
\label{def:idealk}
\rm
Let $K$ denote the ideal of $P$
generated by 
$1-\sum_{i=0}^d x_i$.
\end{definition}

\noindent 
We identify 
$\F\lbrack x_1,  \ldots, x_d\rbrack$
with the $\F$-subalgebra of $P$ generated by
$\lbrace x_i\rbrace_{i=1}^d$.

\begin{lemma}
\label{lem:pk}
The $\F$-vector space $P$ decomposes as
\begin{eqnarray}
P = K+ \F\lbrack x_1,\ldots, x_d\rbrack
\qquad \qquad {\mbox{\rm (direct sum).}}
\label{eq:pk}
\end{eqnarray}
\end{lemma}
\noindent {\it Proof:} 
Let $K_0$ denote the ideal of
$P$ generated by $x_0$. Observe
that the 
$\F$-vector space $P$ decomposes as
\begin{eqnarray}
P = K_0+ \F\lbrack x_1,\ldots, x_d\rbrack
\qquad \qquad {\mbox{\rm (direct sum).}}
\label{eq:pk0}
\end{eqnarray}
Define an $\F$-algebra homomorphism
$\sigma:P\to P$
that sends
$x_0 \mapsto 1-\sum_{i=0}^d x_i$
and fixes $x_j$
for $1 \leq j \leq d$.
The composition of $\sigma $ with itself
is the identity, so $\sigma$ is an isomorphism.
To obtain 
(\ref{eq:pk}),
apply 
$\sigma$ to
each side of (\ref{eq:pk0})
and note that 
$\sigma $ sends
$K_0$ to
$K$ while leaving
$\F\lbrack x_1,\ldots, x_d\rbrack$
invariant.
\hfill $\Box$ \\

\begin{lemma}
\label{lem:pnk}
We have $K\cap P_n=0$ for
$n\geq 0$.
\end{lemma}
 {\it Proof:} 
For notational convenience abbreviate
$y=\sum_{i=0}^d x_i$.
We assume that there exists a nonzero 
$f \in K \cap P_n$ and get a contradiction.
Since $f \in K$ there exists
$h \in P$ such that
$f=(1-y)h$.
Observe that $h\not=0$ since $f\not=0$.
By
(\ref{eq:pdecomp}) there exists an
integer $m\geq 0$ and polynomials
$\lbrace h_i\rbrace_{i=0}^m$  in $P$
such that
$h_i \in P_i$ for $0 \leq i \leq m$
and $h=\sum_{i=0}^m h_i$.
 Without loss we may assume
$h_m\not=0$.
We define some polynomials 
$\lbrace h'_i\rbrace_{i=0}^{m+1}$ as follows:
\begin{eqnarray*}
h'_0=h_0,
\qquad 
h'_i=h_i-yh_{i-1} \;\; (1\leq i \leq m),\qquad
h'_{m+1}=
-yh_m. 
\end{eqnarray*}
Note that 
$h'_i \in P_i$ for
$0 \leq i \leq m+1$, and
$f=
\sum_{i=0}^{m+1}h'_i$.
Observe that $h'_{m+1}\not=0$ since
$P$ is a domain and
each of $y$, $h_m$ is nonzero.
By these comments and since $f\in P_n$ we find
$n=m+1$, $f=h'_{m+1}$, and
$h'_i=0$ for $0 \leq i \leq m$.
Since the $\lbrace h'_i\rbrace_{i=0}^m$ are all zero
we have
$h_0=0$ and
$h_i=yh_{i-1}$ for
$1 \leq i \leq m$. Therefore
$h_i=0$ for $0 \leq i \leq m$.  In particular
$h_m=0$, for a contradiction.
The result follows.
\hfill $\Box$ \\

\begin{lemma}
\label{lem:nuonto}
The following 
are equal:
\begin{enumerate}
\item[\rm (i)]
the image of $K$ under 
$\tilde \nu$;
 \item[\rm (ii)]
the ideal of $\vare^*_0{\tilde T}\vare^*_0$ generated
by $\vare^*_0\Delta \vare^*_0$.
\end{enumerate}
\end{lemma}
{\it Proof:} 
By Corollary
\ref{cor:surtilde}
and
Definition
\ref{def:del},
the image of
$1-\sum_{i=0}^d x_i$ under
$\tilde \nu$ is equal to
$\vare^*_0 \Delta \vare^*_0$.
By Lemma
\ref{lem:tmu}
the map
$\tilde \nu$ is surjective.
The result follows from these comments and Definition
\ref{def:idealk}. 
\hfill $\Box$ \\

\begin{lemma}
\label{lem:diag}
Let
$\pi'$ denote the restriction of
$\pi$
to 
$\vare^*_0{\tilde T}\vare^*_0$.
Let 
$\iota : 
\F\lbrack x_1, \ldots, x_d\rbrack \to P$
denote the inclusion map.
Then the following diagram commutes:

\[
\begin{CD}
P @<\iota<<
         \F\lbrack x_1,\ldots,x_d\rbrack \\
	  @V\tilde \nu VV                     @VV\nu V \\
	\vare^*_0{\tilde T}\vare^*_0 
	  @>>\pi'>
	           e^*_0Te^*_0 
		    \end{CD}
		    \]

\end{lemma}

\noindent {\it Proof:} 
For $1 \leq i \leq d$ we chase $x_i$ around
the diagram.
The image of $x_i$ under
the composition
${\tilde \nu} \circ \iota $ 
is 
$\vare^*_0 \vare_i \vare^*_0$,
and the image of this under
$\pi'$ is
$e^*_0e_ie^*_0$.
The image of
$x_i$ under
$\nu$ is 
$e^*_0e_ie^*_0$.
The result follows.
\hfill $\Box$ \\

\begin{proposition}
\label{conj1}
The following 
are equivalent:
\begin{enumerate}
\item[\rm (i)]
the map $\nu $
is an isomorphism; 
\item[\rm (ii)]
the map ${\tilde \nu} $
is an isomorphism.
\end{enumerate}
\end{proposition}
\noindent {\it Proof:} 
${\rm (i)}\Rightarrow {\rm (ii)}$
The map $\tilde \nu$ is surjective
by
Lemma
\ref{lem:tmu}.
We show
 that $\tilde \nu$ is injective.
By Lemma
\ref{lem:kerintersect} it suffices to 
show that 
$\tilde \nu$ is injective on
$P_n$ for
all integers $n\geq 0$.
Let $n$ be given, and pick any
$f \in P_n$ such that
${\tilde \nu}(f)=0$.
We show $f=0$.
Invoking
Lemma
\ref{lem:pk} we write
$f=k+h$ with $k \in K$
and 
$h \in \F\lbrack x_1,\ldots, x_d\rbrack$.
In the equation
$f=k+h$ 
we apply the composition
$\pi \circ {\tilde \nu}$ to
each term.
The image of
$f$ under 
$\pi \circ {\tilde \nu}$ is zero since
${\tilde \nu}(f)=0$.
The image of
$k$ under 
$\pi \circ {\tilde \nu}$ is zero
by Lemma
\ref{lem:nuonto} and
Proposition
\ref{thm:idealinter}(i),(iv).
The image of $h$ under
$\pi \circ {\tilde \nu}$ is 
$\nu(h)$ by 
Lemma
\ref{lem:diag}.
 By these comments
$\nu(h)=0$.
We assume
$\nu$ is an isomomorphism
so $h=0$.
Therefore $f=k \in K$.
We have $f \in K$ and
$f\in P_n$,
so $f=0$ in view of
Lemma
\ref{lem:pnk}.
\\
${\rm (ii)}\Rightarrow {\rm (i)}$
The map
$\nu$ is surjective by
Corollary
\ref{cor:nusurj}.
We show 
that $\nu$ is injective.
Suppose we are given 
$h \in 
\F \lbrack x_1,\ldots, x_d\rbrack$
such that $\nu(h)=0$.
We show $h=0$.
By Lemma
\ref{lem:diag} the composition
$\pi \circ {\tilde \nu}$ sends
$h \mapsto 0$.
Therefore
${\tilde \nu}(h)$ is in the kernel of
$\pi$.
By this and
Proposition
\ref{thm:idealinter}(i),(iv)
we see that
${\tilde \nu}(h)$ is in the ideal of
 $\vare^*_0{\tilde T} \vare^*_0$
 generated by
 $\vare^*_0{\Delta} \vare^*_0$.
Now
$h \in K$ by
 Lemma
\ref{lem:nuonto}
and since $\tilde \nu$ is an isomorphism.
We have $h \in K$ and
$h \in \F \lbrack x_1,\ldots, x_d\rbrack$,
so $h=0$ in view of
Lemma
\ref{lem:pk}.
\hfill $\Box$ \\

\begin{corollary}
\label{cor:iff}
The following 
are equivalent:
\begin{enumerate}
\item[\rm (i)]
the map $\mu $
from
Corollary \ref{cor:18a}
is an isomorphism; 
\item[\rm (ii)]
the map ${\tilde \nu} $
from Corollary
\ref{cor:surtilde}
is an isomorphism.
\end{enumerate}
\end{corollary}
\noindent {\it Proof:} 
Combine
Proposition
\ref{lem:isodual}
and
Proposition
\ref{conj1}.
\hfill $\Box$ \\

\section{The proof of Theorem \ref{conj:main} and
Theorem 
\ref{conj:mainp}
}

\noindent In this section we prove
 Theorem \ref{conj:main} and
Theorem 
\ref{conj:mainp}.

\medskip
\noindent Throughout this section fix an integer $d \geq 0$.
Recall the sets 
${\rm Feas}(d,\F)$ 
from Definition
\ref{def:feaset}
and 
${\rm Rac}(d,\F)$
from Definition
\ref{def:RAC}.

\begin{definition}
\label{def:good}
\rm
Pick any sequence $ p \in 
{\rm Feas}(d,\F)$  and consider the ordered pair
$(p, \F)$. This pair is said to be {\it confirmed}
whenever Theorem
\ref{conj:mainp} is true for that
$p$ and $\F$.  
\end{definition}

\begin{lemma}
\label{lem:step5}
{\rm 
\cite[Theorem~12.1]{nom:mu}}
Assume $d \leq  5$. Then $(p,\F)$ is confirmed for
all 
$p \in{\rm Feas}(d,\F)$.
\end{lemma}

\begin{lemma}
\label{lem:qrstep}
{\rm 
\cite[Theorem~5.3]{NT:muqrac}}
The pair 
 $(p,\F)$ is confirmed
 for all 
$p \in{\rm Rac}(d,\F)$.
\end{lemma}

\begin{lemma}
\label{lem:fieldchange}
{\rm 
\cite[Theorem~5.2]{NT:muqrac}}
Pick any
$p \in{\rm Feas}(d,\F)$. If 
there exists a field extension
$\Knew$ 
of $\F$ such that 
$(p,\Knew)$  is confirmed, then
$(p,\F)$  is confirmed. 
\end{lemma}

\noindent {\it Proof of Theorem
\ref{conj:mainp}:}
We will confirm the pair $(p,\F)$ in the sense
of 
Definition
\ref{def:good}.
We may assume $d\geq 3$; otherwise 
 $(p,\F)$ is confirmed by
 Lemma
\ref{lem:step5}.
Abbreviate $\Knew={\overline \F} $ for the algebraic
closure of $\F$, and note that $\Knew$ is infinite.
By Lemma
\ref{lem:qrstep} the pair
$(p',\Knew)$ is confirmed for all
sequences
$p' \in{\rm Rac}(d,\Knew)$.
Now by 
Proposition
\ref{prop:equiv} and
Corollary
\ref{cor:iff},
the type $\lbrack n \rbrack$ is
$(p',\Knew)$-direct
for all integers $n\geq 0$ and
all sequences $p' \in{\rm Rac}(d,\Knew)$.
Now by
Proposition
\ref{prop:sometimes} and since
$\Knew $ is infinite,
the type $\lbrack n \rbrack$ is
$(p',\Knew)$-direct
for all integers $n\geq 0$ and all
sequences $p' \in{\rm Feas}(d,\Knew)$.
In particular
the type $\lbrack n \rbrack$ is
$(p,\Knew)$-direct
for all integers $n\geq 0$.
Now by 
Proposition
\ref{prop:equiv} and
Corollary
\ref{cor:iff},
the pair $(p,\Knew)$ is confirmed.
Now by
Lemma
\ref{lem:fieldchange}
the pair $(p,\F)$ is confirmed.
\hfill $\Box$ \\

\noindent {\it Proof of Theorem
 \ref{conj:main}:  
}
Immediate from
Theorem
\ref{conj:mainp}
and 
\cite[Theorem~10.1]{nom:mu}.
\hfill $\Box$ \\

\section{Comments}

\noindent In the previous section we proved 
 Theorem \ref{conj:main} and
Theorem 
\ref{conj:mainp}.
In this section
we list some related results that
might be of independent interest. We also mention a
conjecture.

\medskip
\noindent The following is a corollary
to Theorem
\ref{conj:main}.

\begin{corollary}
 \label{conj:mainac}
 Assume the field $\F$ is algebraically closed.
Let $d$ denote a nonnegative integer and  let
\begin{equation}         \label{eq:parrayac}
 (\{\theta_i\}_{i=0}^d; \{\theta^*_i\}_{i=0}^d; \{\zeta_i\}_{i=0}^d)
\end{equation}
denote a sequence of scalars taken from $\K$.
Then there exists a TD system $\Phi$ over $\K$
with parameter array \eqref{eq:parrayac} if and only if
{\rm (i)--(iii)} hold below.
\begin{itemize}
\item[\rm (i)]
$\theta_i \neq \theta_j$, 
$\theta^*_i \neq \theta^*_j$ if $i \neq j$ $(0 \leq i,j \leq d)$.
\item[\rm (ii)]
The expressions
\begin{equation*} 
\label{eq:betaplusoneac} 
\frac{\theta_{i-2}-\theta_{i+1}}{\theta_{i-1}-\theta_i},  \qquad\qquad
  \frac{\theta^*_{i-2}-\theta^*_{i+1}}{\theta^*_{i-1}-\theta^*_i}
\end{equation*}
are equal and independent of $i$ for $2 \leq i \leq d-1$.
\item[\rm (iii)]
$\zeta_0=1$, $\zeta_d \neq 0$, and
\begin{equation*}   
        \label{eq:ineqac}
0 \neq \sum_{i=0}^d \eta_{d-i}(\theta_0)\eta^*_{d-i}(\theta^*_0) \zeta_i.
\end{equation*}
\end{itemize}
Suppose {\rm (i)--(iii)} hold. Then $\Phi$ is unique up to isomorphism of
TD systems.
\end{corollary}
\noindent {\it Proof:}
By
Theorem
\ref{conj:main} and since
every tridiagonal system over 
an algebraically closed field 
 is sharp
\cite[Theorem~1.3]{nomstructure}.
\hfill $\Box$ \\

\begin{theorem}
Fix an integer $d\geq 0$ and a 
  sequence $p \in
{\rm Feas}(d,\F)$.
Let the algebra
$T=T(p,\F)$ 
be as in Definition
\ref{def:twitha}. Then the corresponding map
$\nu: \F\lbrack x_1,\ldots, x_d\rbrack
\to e^*_0Te^*_0$ from
Definition
\ref{def:nu}
is an isomorphism.
\end{theorem}
\noindent {\it Proof:}
Combine
Theorem
\ref{conj:mainp}
and
Proposition
\ref{lem:isodual}.
\hfill $\Box$ \\

\begin{theorem}
\label{thm:dirnu}
Fix an integer $d\geq 0$ and a 
  sequence $p \in
{\rm Feas}(d,\F)$.
Let the algebra
${\tilde T}={\tilde T}(p,\F)$ 
be as in Definition
\ref{def:tt}.
Then the corresponding map
${\tilde \nu}: P
\to \vare^*_0{\tilde T} \vare^*_0$ from
Corollary
\ref{cor:surtilde}
is an isomorphism.
\end{theorem}
\noindent {\it Proof:}
Combine
Theorem
\ref{conj:mainp}
and
Corollary
\ref{cor:iff}.
\hfill $\Box$ \\

\begin{theorem}
Fix an integer $d\geq 0$ and a 
  sequence $p \in
{\rm Feas}(d,\F)$.
Then for all integers $n\geq 0$ the
type $\lbrack n \rbrack$ is
$(p,\F)$-direct in the sense of
Definition \ref{def:fav}.
\end{theorem}
\noindent {\it Proof:}
Combine
Proposition
\ref{prop:equiv}
and 
Theorem
\ref{thm:dirnu}.
\hfill $\Box$ \\

\begin{theorem}
\label{thm:spanzz}
Fix an integer $d\geq 0$ and a 
  sequence $p \in
{\rm Feas}(d,\F)$.
Let the algebra
${\tilde T}={\tilde T}(p,\F)$ 
be as in Definition
\ref{def:tt}.
Then the $\F$-vector space $\tilde T$ is spanned by its
zigzag words.
\end{theorem}
\noindent {\it Proof:}
Combine
Proposition
\ref{lem:tildehc}
and Proposition
\ref{thm:zz}.
\hfill $\Box$ \\

\noindent Below Lemma
\ref{cor:zplusr}
we conjectured that the sum
(\ref{eq:dirsum}) is always direct.
In the context of $\tilde T$ this conjecture
can be expressed as follows.

\begin{conjecture}
\label{conj:final}
\rm
Fix an integer $d\geq 0$ and a 
  sequence $p \in
{\rm Feas}(d,\F)$.
Let the algebra
${\tilde T}={\tilde T}(p,\F)$ 
be as in Definition
\ref{def:tt}.
Then the $\F$-vector space $\tilde T$ has a basis consisting of
its 
zigzag words.
\end{conjecture}

\bigskip

\noindent Tatsuro Ito \hfil\break
\noindent Division of Mathematical and Physical Sciences \hfil\break
\noindent Graduate School of Natural Science and Technology\hfil\break
\noindent Kanazawa University \hfil\break
\noindent Kakuma-machi,  Kanazawa 920-1192, Japan \hfil\break
\noindent email:  {\tt tatsuro@kenroku.kanazawa-u.ac.jp} \hfil\break

\bigskip

\noindent Kazumasa Nomura
\hfil\break
\noindent
626-1-109, Awano, Kamagaya-shi
\hfil\break
\noindent
Chiba, 273-0132 Japan
\hfil\break
\noindent
email: {\tt knomura@pop11.odn.ne.jp}
\hfil\break

\bigskip

\noindent Paul Terwilliger \hfil\break
\noindent Department of Mathematics \hfil\break
\noindent University of Wisconsin \hfil\break
\noindent 480 Lincoln Drive \hfil\break
\noindent Madison, WI 53706-1388 USA \hfil\break
\noindent email: {\tt terwilli@math.wisc.edu }\hfil\break
\end{document}